# NONPARAMETRIC QUASI-MAXIMUM LIKELIHOOD ESTIMATION FOR GAUSSIAN LOCALLY STATIONARY PROCESSES[1]

By Rainer Dahlhaus and Wolfgang Polonik

*Universität Heidelberg and University of California, Davis*

This paper deals with nonparametric maximum likelihood estimation for Gaussian locally stationary processes. Our nonparametric MLE is constructed by minimizing a frequency domain likelihood over a class of functions. The asymptotic behavior of the resulting estimator is studied. The results depend on the richness of the class of functions. Both sieve estimation and global estimation are considered.

Our results apply, in particular, to estimation under shape constraints. As an example, autoregressive model fitting with a monotonic variance function is discussed in detail, including algorithmic considerations.

A key technical tool is the time-varying empirical spectral process indexed by functions. For this process, a Bernstein-type exponential inequality and a central limit theorem are derived. These results for empirical spectral processes are of independent interest.

**1. Introduction.** Nonstationary time series whose behavior is locally close to the behavior of a stationary process can often be successfully described by models with time-varying parameters, that is, by models characterized by parameter curves. A simple example is the time-varying AR(1) model $X_t + \alpha_t X_{t-1} = \sigma_t \varepsilon_t$, $t \in \mathbf{Z}$, where $\alpha_t$ and $\sigma_t$ vary over time. If the process is observed at times $t = 1, \ldots, n$, the problem of estimation of $\alpha_t$ and $\sigma_t$ may be formulated as the estimation of the curves $\alpha(\cdot)$ and $\sigma(\cdot)$ with $\alpha(t/n) = \alpha_t$, $\sigma(t/n) = \sigma_t$ in an adequately rescaled model. To study such problems in a

Received November 2003; revised September 2005.
[1]Supported in part by the Deutsche Forschungsgemeinschaft under DA 187/12-2 (R. D.) and by NSF Grants 0103606 and 0406431 (W. P.).

more general framework, Dahlhaus [9] introduced the class of locally stationary processes having a time-varying spectral representation or, alternatively, an infinite time-varying moving average representation. In this paper, we present a methodology for nonparametric ML-estimation of time-varying spectral densities of Gaussian locally stationary processes. Results for parameter functions like $\alpha(\cdot)$ or $\sigma(\cdot)$ then follow from the former results for spectral densities. The time-varying AR(1)-process from above will serve as a simple example for our general results.

Guo et al. [16] consider an approach for nonparametric estimation of the time-varying spectral density using both a penalized least squares and a penalized likelihood approach. For nonparametric estimation of curves such as $\alpha(\cdot)$ and $\sigma(\cdot)$ in the above example, different approaches have been considered. One idea is to utilize a stationary method on overlapping small segments of the time series (e.g., a Yule–Walker or least squares estimate) where the resulting estimate is regarded as the estimate of the curve at the midpoint of the interval. More generally, one can consider kernel estimates [11] or local linear fits, as in [17]. Other methods are based on wavelets, as in [12] and in [13].

In contrast to these local methods, we here consider a global method by fitting parameter curves from a given class of functions. Such a method is of particular interest when shape restrictions are known, as, for instance, in case of earthquake data or chirp signals where some of the parameter functions are known to be monotonic or unimodal (cf. Section 3). We fit such curves by maximizing an appropriate likelihood function over a class of suitable candidate functions. By choosing the class of functions in an adequate way, different estimates can be obtained. We consider both sieve estimates and global estimates in function spaces. The likelihood used is a minimum distance functional between spectral densities in the "time-frequency" domain, meaning that the spectral densities are functions of both time and frequency. The likelihood considered here can be regarded as a generalization of the classical Whittle likelihood [24] to locally stationary processes.

The basic technical tool for deriving rates of convergence of the nonparametric maximum likelihood estimator is an exponential inequality for the time-varying empirical spectral process of a locally stationary process (cf. [14]).

Non- and semiparametric inference has received a lot of attention during the last decade. A general approach uses minimum contrast estimation, where some contrast functional is minimized over an infinite-dimensional parameter space, including maximum likelihood estimation, $M$-estimation, least squares estimates in nonparametric regression (e.g., [2, 3, 7, 21, 22]). The theory for all of these approaches is based on the behavior of some kind of empirical process whose analysis crucially depends on exponential inequalities (or concentration inequalities) together with measures of complexity of the parameter space, such as metric entropy conditions or VC indices. The theory often leads to (almost) optimal rates of convergence for the estimates.

It turns out that by using our approach to nonparametric ML-estimation for locally stationary processes, it is possible to follow some of the main steps of the approaches mentioned above. However, the statistical problem, the likelihood under consideration, the underlying empirical process and, hence, the technical details, are quite different. For instance, our contrast functional turns out to be equivalent to an $L_2$-distance in the time-frequency domain (instead of the Hellinger distance, as in the case of van de Geer). Further, we do not exploit metric entropy with bracketing, since the time-varying empirical spectral process is not monotone in its argument. This, in fact, led us to also consider sieve estimation. In addition, there is, of course, the complex dependence structure for locally stationary processes which, for example, enters when proving exponential inequalities for the increments of the empirical spectral process.

In Section 2, we describe the estimation problem and the construction of the likelihood and we present the main results on rate of convergence of our nonparametric likelihood estimates. In Section 3, the estimation of a monotonic variance function in a time-varying AR-model is studied, including explicit algorithms involving isotonic regression. In Section 4, we prove a Bernstein-type exponential inequality for the function-indexed empirical spectral process for Gaussian locally stationary processes. This exponential inequality is used to derive maximal inequalities and a functional limit theorem. All proofs are shifted without further reference to Section 5. The Appendix contains some auxiliary results.

## 2. Basic ideas and the main result.

2.1. *Locally stationary processes.* In this paper, we assume that the observed process is Gaussian and *locally stationary*. Locally stationary processes were introduced in [9] by using a time-varying spectral representation. In contrast to this, we use a time-varying MA($\infty$)-representation and formulate the assumptions in the time domain. As in nonparametric regression, we rescale the functions in time to the unit interval in order to achieve a meaningful asymptotic theory. The setup is more general than, for example, in [9] since we allow for jumps in the parameter curves by assuming bounded variation instead of continuity in the time direction.

Let

$$V(g) = \sup\left\{\sum_{k=1}^{m} |g(x_k) - g(x_{k-1})| : 0 \le x_0 < \cdots < x_m \le 1,\ m \in \mathbf{N}\right\}$$

be the total variation of a function $g$ on $[0,1]$, and for some $\kappa > 0$, let

$$\ell(j) := \begin{cases} 1, & |j| \leq 1, \\ |j|\log^{1+\kappa}|j|, & |j| > 1. \end{cases}$$

DEFINITION 2.1 (Locally stationary processes). The sequence $X_{t,n}, t = 1, \ldots, n$, is a *locally stationary process* if it has the representation

(1) $$X_{t,n} = \sum_{j=-\infty}^{\infty} a_{t,n}(j)\varepsilon_{t-j},$$

where the $\varepsilon_t$ are identically distributed with $E\varepsilon_t \equiv 0$, $E\varepsilon_s\varepsilon_t = 0$ for $s \neq t$, $E\varepsilon_t^2 \equiv 1$ and where the following conditions hold:

(2) $$\sup_t |a_{t,n}(j)| \leq \frac{K}{\ell(j)} \qquad \text{(with $K$ not depending on $n$)}$$

and there exist functions $a(\cdot, j):(0,1] \to \mathbf{R}$ with

(3) $$\sup_u |a(u,j)| \leq \frac{K}{\ell(j)},$$

(4) $$\sup_j \sum_{t=1}^n \left| a_{t,n}(j) - a\left(\frac{t}{n},j\right)\right| \leq K,$$

(5) $$V(a(\cdot,j)) \leq \frac{K}{\ell(j)}.$$

If the process $X_{t,n}$ is Gaussian (as in this paper), it can be shown that the $\varepsilon_t$ also have to be Gaussian.

The above conditions are discussed in [14]. A simple example of a process $X_{t,n}$ which fulfills the above assumptions is $X_{t,n} = \phi(\frac{t}{n})Y_t$, where $Y_t = \Sigma_j a(j)\varepsilon_{t-j}$ is stationary with $|a(j)| \leq K/\ell(j)$ and $\phi$ is of bounded variation. In [14], Theorem 2.3, we have shown that time-varying ARMA (tvARMA) models whose coefficient functions are of bounded variation are locally stationary in the above sense. In particular, it follows from this result that the system of difference equations

(6) $$X_{t,n} + \sum_{j=1}^p \alpha_j\left(\frac{t}{n}\right)X_{t-j,n} = \sigma\left(\frac{t}{n}\right)\varepsilon_t,$$

where $\varepsilon_t$ are i.i.d. with $E\varepsilon_t = 0$ and $E|\varepsilon_t| < \infty$, all $\alpha_j(\cdot)$ as well as $\sigma^2(\cdot)$ are of bounded variation, $1 + \sum_{j=1}^p \alpha_j(u)z^j \neq 0$ for all $u$ and all $z$ such that $0 < |z| \leq 1 + \delta$ for some $\delta > 0$, has a locally stationary solution which is called *tvAR process*.

DEFINITION 2.2 (Time-varying spectral density and covariance). Let $X_{t,n}$ be a locally stationary process. The function

$$f(u,\lambda) := \frac{1}{2\pi}|A(u,\lambda)|^2$$

with

$$A(u,\lambda) := \sum_{j=-\infty}^{\infty} a(u,j)\exp(-i\lambda j)$$

is the *time-varying spectral density*, and

(7) $$c(u,k) := \int_{-\pi}^{\pi} f(u,\lambda)\exp(i\lambda k)\,d\lambda = \sum_{j=-\infty}^{\infty} a(u,k+j)a(u,j)$$

is the *time-varying covariance* of lag $k$ at rescaled time $u$.

For instance, the time-varying spectral density of a tvAR($p$) process is given by

(8) $$f(u,\lambda) = \frac{\sigma^2(u)}{2\pi}\left|1 + \sum_{j=1}^{p} \alpha_j(u)\exp(i\lambda j)\right|^{-2}.$$

2.2. *The estimator.* Our nonparametric estimator of the time-varying spectral density of a locally stationary process will be defined as a minimum contrast estimator in the time-frequency domain, that is, we minimize a contrast functional between a nonparametric estimate of the time-varying spectral density $\hat{f}(u,\lambda)$ over a class of candidate spectral density functions. Two different scenarios are considered: (i) Sieve estimation, where the classes of candidate spectral densities $\mathcal{F}_n$ depend on $n$, are "finite-dimensional" and approximate, as $n$ gets large, a (large) target class $\mathcal{F}$ and (ii) global estimation, where the contrast functional is minimized over an "infinite-dimensional" target class $\mathcal{F}$ directly, or, formally, $\mathcal{F}_n = \mathcal{F}$ for all $n$. The nonparametric (sieve) maximum likelihood estimate for $f$ is defined by

$$\widehat{f}_n = \underset{g \in \mathcal{F}_n}{\operatorname{argmin}}\, \mathcal{L}_n(g),$$

where our contrast functional is

(9) $$\mathcal{L}_n(g) = \frac{1}{n}\sum_{t=1}^{n} \frac{1}{4\pi}\int_{-\pi}^{\pi} \left\{\log g\left(\frac{t}{n},\lambda\right) + \frac{J_n(\frac{t}{n},\lambda)}{g(\frac{t}{n},\lambda)}\right\} d\lambda.$$

Here, $J_n$ denotes a nonparametric estimate of the time-varying spectral density, defined as

(10) $$J_n\left(\frac{t}{n},\lambda\right) = \frac{1}{2\pi}\sum_{k:\,1\leq[t+1/2\pm k/2]\leq n} X_{[t+1/2+k/2],n} X_{[t+1/2-k/2],n}\exp(-i\lambda k).$$

This estimate is called pre-periodogram [19]. It can be regarded as a preliminary estimate of $f(\frac{t}{n}, \lambda)$; however, in order to become consistent, it has to be smoothed in the time and frequency directions. If we choose $g(u, \lambda) = \tilde{g}(\lambda)$, that is, we model the underlying time series as stationary, then $\mathcal{L}_n(g)$ is identical to the classical Whittle likelihood. This follows since the classical periodogram is a time average of the pre-periodogram.

Below, we prove convergence of $\hat{f}_n$ to

(11) $$f_{\mathcal{F}} = \operatorname*{argmin}_{g \in \mathcal{F}} \mathcal{L}(g),$$

where

(12) $$\mathcal{L}(g) = \int_0^1 \frac{1}{4\pi} \int_{-\pi}^{\pi} \left\{ \log g(u, \lambda) + \frac{f(u, \lambda)}{g(u, \lambda)} \right\} d\lambda \, du.$$

This is, up to a constant, the asymptotic Kullback–Leibler information divergence between two Gaussian locally stationary processes with mean zero and time-varying spectra $f(u, \lambda)$ and $g(u, \lambda)$ (cf. [8], Theorem 3.4). Since $\mathcal{L}(g) \geq \int_0^1 \frac{1}{4\pi} \int_{-\pi}^{\pi} \{\log f(u, \lambda) + 1\} d\lambda \, du$, we have

$$f_{\mathcal{F}} = f \iff f \in \mathcal{F},$$

provided the minimizer in (11) is unique (a.s., uniqueness of the minimizer follows in the case $f \in \mathcal{F}$ from the inequality $\log x < x - 1 \ \forall x \neq 1$).

We now give three examples of possible model classes $\mathcal{F}$. In these examples and in all of what follows, candidate spectral densities are denoted by $g(u, \lambda)$. The true spectral density is always denoted by $f(u, \lambda)$.

EXAMPLE 2.3 (Model classes for the time-varying spectrum). (a) The locally stationary process is parameterized by curves $\theta(u) = (\theta_1(u), \ldots, \theta_d(u))' \in \Theta$ and $\mathcal{F}$ consists of all spectral densities of the form $g_\theta(u, \lambda) = w(\theta(u), \lambda)$ for some fixed function $w$ as, for instance, in the case of tvAR models discussed above.

(b) The process is stationary, that is, $g(u, \lambda) = \tilde{g}(\lambda)$, and the spectral density $\tilde{g}(\lambda)$ is the curve to be estimated nonparametrically.

(c) Both the behavior in time and in frequency are modeled nonparametrically. An example is the amplitude-modulated process $X_{t,n} = \varphi(\frac{t}{n}) Y_t$ where $Y_t$ is stationary. In this case, $g(u, \lambda) = g_1(u) g_2(\lambda)$, where $g_1(u) = \varphi(u)^2$ and $g_2(\lambda)$ is the spectral density of $Y_t$.

In the above examples, the curves $\theta_j(u)$, $\tilde{g}(\lambda)$, $g_1(u)$ and $g_2(\lambda)$ are assumed to lie in 'smoothness' classes, like Sobolev classes or classes defined through shape restrictions (see Section 3 below).

2.3. *Asymptotic properties of the NPMLE.* We now motivate and formulate rates of convergence for the NPMLE. It turns out that sieve estimation leads to rates of convergence for the NPMLE which, in an i.i.d. density estimation setup, are known to be (almost) optimal. In the case of global estimation, the obtained rates are slower. Whether the same rates as for sieve estimation can be obtained for global estimation is an open question.

We start with some elementary calculations which demonstrate the structure of the problem, the importance of the empirical spectral process and the fact that the $L_2$-norm of the inverse spectral densities is a natural norm for studying the convergence of the NPMLE.

First, we define the empirical spectral process by

$$(13) \qquad E_n(\phi) = \sqrt{n}(F_n(\phi) - F(\phi)),$$

where

$$(14) \qquad F(\phi) = \int_0^1 \int_{-\pi}^{\pi} \phi(u,\lambda) f(u,\lambda)\, d\lambda\, du$$

and

$$(15) \qquad F_n(\phi) = \frac{1}{n} \sum_{t=1}^{n} \int_{-\pi}^{\pi} \phi\left(\frac{t}{n}, \lambda\right) J_n\left(\frac{t}{n}, \lambda\right) d\lambda.$$

In the following motivation, we only consider the case $\mathcal{F}_n = \mathcal{F}$ for all $n$. The case of sieve estimation is similar in nature [see proof of Theorem 2.6, part (I)]. By definition of $\widehat{f}_n$ and $f_\mathcal{F}$, we have

$$(16) \qquad \mathcal{L}_n(\widehat{f}_n) \leq \mathcal{L}_n(f_\mathcal{F})$$

and, similarly,

$$(17) \qquad \mathcal{L}(f_\mathcal{F}) \leq \mathcal{L}(\widehat{f}_n).$$

Combining (16) and (17), we obtain the basic inequalities

$$(18) \qquad \begin{aligned} 0 \leq \mathcal{L}(\widehat{f}_n) - \mathcal{L}(f_\mathcal{F}) &\leq (\mathcal{L}_n - \mathcal{L})(f_\mathcal{F}) - (\mathcal{L}_n - \mathcal{L})(\widehat{f}_n) \\ &= \frac{1}{4\pi} \frac{1}{\sqrt{n}} E_n\left(\frac{1}{f_\mathcal{F}} - \frac{1}{\widehat{f}_n}\right) + R_{\log}(f_\mathcal{F}) - R_{\log}(\widehat{f}_n), \end{aligned}$$

where

$$(19) \qquad R_{\log}(g) := \frac{1}{4\pi} \int_{-\pi}^{\pi} \left[ \frac{1}{n} \sum_{t=1}^{n} \log g\left(\frac{t}{n}, \lambda\right) - \int_0^1 \log g(u,\lambda)\, du \right] d\lambda.$$

Hence, if $\sup_{g \in \mathcal{F}} |R_{\log}(g)|$ is small, the convergence of $\mathcal{L}(\widehat{f}_n) - \mathcal{L}(f_\mathcal{F})$ can be controlled by the empirical spectral process whose properties will be investigated in Section 4, leading to the subsequent convergence results. Note that in the correctly specified case where $f_\mathcal{F} = f$,

$$(20) \quad \mathcal{L}(\widehat{f}_n) - \mathcal{L}(f_\mathcal{F}) = \frac{1}{4\pi} \int_0^1 \int_{-\pi}^{\pi} \left\{ \log \frac{\widehat{f}_n(u,\lambda)}{f(u,\lambda)} + \frac{f(u,\lambda)}{\widehat{f}_n(u,\lambda)} - 1 \right\} d\lambda \, du,$$

which equals the Kullback–Leibler information divergence between two Gaussian locally stationary processes (cf. [8]). Under certain assumptions, the equivalence of the above information divergence to $\rho_2(1/\widehat{f}_n, 1/f)^2$ is shown below (Lemma 5.1), where $\rho_2(\phi, \psi) = \rho_2(\phi - \psi)$ with

$$(21) \qquad \rho_2(\phi) = \left( \int_0^1 \int_{-\pi}^{\pi} |\phi(u,\lambda)|^2 \, d\lambda \, du \right)^{1/2}.$$

Hence, properties of the empirical spectral process lead, via (18), to the convergence results for $\rho_2^2(1/\widehat{f}_n, 1/f)$ stated in Theorem 2.6. The above discussion shows that these results are also convergence results for the Kullback–Leibler information divergence. More generally, we allow for misspecification in our results, which means that we do not require $f_\mathcal{F} = f$. In this case, additional convexity arguments come into play (cf. Lemma 5.1). In order to formulate the assumptions on the class $\mathcal{F}$, we need to introduce further notation. With

$$(22) \qquad \hat{\phi}(u,j) := \int_{-\pi}^{\pi} \phi(u,\lambda) \exp(i\lambda j) \, d\lambda,$$

let

$$(23) \quad \begin{aligned} \rho_\infty(\phi) &:= \sum_{j=-\infty}^{\infty} \sup_u |\hat{\phi}(u,j)|, \quad \tilde{v}(\phi) := \sup_j V(\hat{\phi}(\cdot,j)) \quad \text{and} \\ v_\Sigma(\phi) &:= \sum_{j=-\infty}^{\infty} V(\hat{\phi}(\cdot,j)). \end{aligned}$$

Furthermore, let

$$\mathcal{F}^* = \left\{ \frac{1}{g}; g \in \mathcal{F} \right\}.$$

Since the empirical process has $1/g$ as its argument, it is more natural to use the class $\mathcal{F}^*$ instead of $\mathcal{F}$ in most of the assumptions. For our sieve estimate, we also need a sequence of approximating classes denoted by $\mathcal{F}_n$. The corresponding classes of inverse functions are denoted by $\mathcal{F}_n^*$. In the results on global estimation, $\mathcal{F}_n^* \equiv \mathcal{F}^*$.

ASSUMPTION 2.4. (a) The classes $\mathcal{F}_n$ are such that $\widehat{f}_n$ exists for all $n$, and $\mathcal{F}$ is such that $f_\mathcal{F}$ exists, is unique and $0 < f_\mathcal{F} < \infty$.

(b) For any $\phi \in \mathcal{F}^*$, there exists a sequence $\pi_n(\phi) \in \mathcal{F}_n^*$ such that $\rho_2(\phi, \pi_n(\phi)) \to 0$ as $n \to \infty$.

(c) There exist $0 < M_* \leq 1 \leq M^* < \infty$ with $M_* \leq |\phi(u, \lambda)| < M^*$ for all $u, \lambda$ and $\phi \in \mathcal{F}_n^*$. Furthermore, $\sup_{\phi \in \mathcal{F}_n^*} \rho_\infty(\phi) \leq \rho_\infty < \infty$, $\sup_{\phi \in \mathcal{F}_n^*} \tilde{v}(\phi) < \tilde{v} < \infty$ and $\sup_{\phi \in \mathcal{F}_n^*} v_\Sigma(\phi) \leq v_\Sigma < \infty$. All constants may depend on $n$.

The bounds in (c) are not very restrictive. For instance, for tvAR processes, only finitely many $\widehat{\phi}(u, j)$ are different from zero; see Example 2.7 below. $M_*$ (the uniform upper bound on the model spectra) is only needed for bounding $R_{\log}(g)$ in Lemma A.2. This bound can be avoided by a condition on the variation of $\int_{-\pi}^{\pi} \log g(u, \lambda) \, d\lambda$ which in some cases already follows from other assumptions; see Example 2.7 and Section 3.1 below. In that case, the constant $M_*$ in Theorem 2.6 can be replaced by 1. We mention the following elementary relationships:

$$\sup_{u,\lambda} |\phi(u,\lambda)| \leq \frac{1}{2\pi} \rho_\infty(\phi), \qquad \tilde{v}(\phi) \leq v_\Sigma(\phi),$$

$$\tilde{v}(\phi) \leq \int V(\phi(\cdot, \lambda)) \, d\lambda, \qquad \rho_2(\phi) \leq \frac{1}{\sqrt{2\pi}} \rho_\infty(\phi).$$

Our results for the NPMLE are derived under conditions on the richness of the model class $\mathcal{F}^*$, as measured by metric entropy. For each $\epsilon > 0$, the *covering number* of a class of functions $\Phi$ with respect to the metric $\rho_2$ is defined as

(24) $$N(\epsilon, \Phi, \rho_2) = \inf\{n \geq 1 : \exists \phi_1, \ldots, \phi_n \in \Phi \text{ such that } \forall \phi \in \Phi \, \exists 1 \leq i \leq n \text{ with } \rho_2(\phi, \phi_i) < \epsilon\}.$$

The quantity $H(\epsilon, \Phi, \rho_2) = \log N(\epsilon, \Phi, \rho_2)$ is called the *metric entropy* of $\Phi$ with respect to $\rho_2$. For technical reasons, we assume that $H(\epsilon, \Phi, \rho_2) \leq \tilde{H}_\Phi(\epsilon)$ with $\tilde{H}_\Phi(\cdot)$ continuous and monotonically decreasing. This assumption is known to be satisfied for many function classes (see Example 2.7). A crucial quantity is the *covering integral*

(25) $$\int_\epsilon^\delta \tilde{H}_\Phi(u) \, du.$$

In contrast to (25), the standard covering integral is defined to be the integral over the *square root* of the metric entropy. Here, we have to use this larger covering integral which leads to slower rates of convergence as compared to nonparametric ML-estimation based on i.i.d. data (cf. the discussion in Section 2.4).

REMARK 2.5 (Measurability). We will not discuss measurability here. All random quantities considered are assumed to be measurable. In the case where $\mathcal{F}_n^*$ is nonseparable, this measurability assumption may be an additional restriction.

THEOREM 2.6 (Rates of convergence). *Let $X_{t,n}$ be a Gaussian locally stationary process. Let $\mathcal{F}_n$, $\mathcal{F}$ be classes of functions satisfying Assumption* 2.4.

Part I (Sieve estimation) *Suppose that there exist constants $A > 0$, $k_n \geq 1$ with $\log N(\eta, \mathcal{F}_n^*, \rho_2) \leq A k_n \log(n/\eta)$ for all $\eta > 0$. Let $c_n = \max\{\rho_\infty, v_\Sigma, (M^*)^2\}$. If $f \in \mathcal{F}$, we have, with $a_n = \inf_{g \in \mathcal{F}_n} \rho_2(\frac{1}{f_\mathcal{F}}, \frac{1}{g})$, that*

$$\rho_2\left(\frac{1}{\widehat{f}_n}, \frac{1}{f_\mathcal{F}}\right) = O_P(\delta_n),$$

*where $\delta_n$ satisfies*

(26) $$\delta_n = \frac{M^*}{M_*} \max\left(\sqrt{\frac{c_n k_n \log n}{n}}, a_n\right).$$

*If $f \notin \mathcal{F}$, the same result holds, with $a_n$ replaced by $b_n = \rho_2(\frac{1}{f_\mathcal{F}}, \frac{1}{f_{\mathcal{F}_n}})$, provided that all $\mathcal{F}_n^*$ are convex.*

Part II (Global estimation) *Let $\mathcal{F}_n = \mathcal{F}$. Assume either $f \in \mathcal{F}$ or $\mathcal{F}^*$ to be convex. Further, assume that there exist $0 < \gamma < 2$ and $0 < A < \infty$ such that for all $\eta > 0$,*

(27) $$\tilde{H}_{\mathcal{F}^*}(\eta) \leq A\eta^{-\gamma}.$$

*Then*

$$\rho_2\left(\frac{1}{\widehat{f}_n}, \frac{1}{f_\mathcal{F}}\right) = O_P(\delta_n),$$

*where*

$$\delta_n = \begin{cases} n^{-\frac{1}{2(\gamma+1)}} & \text{for } 0 < \gamma < 1, \\ n^{-\frac{2-\gamma}{4\gamma}} (\log n)^{\frac{\gamma-1}{2\gamma}} & \text{for } 1 \leq \gamma < 2. \end{cases}$$

*Remark.* In Part I, the nonrandom term $a_n$ is smaller than $b_n$. Furthermore, (an upper bound of) $a_n$ may be easier to calculate.

EXAMPLE 2.7. The above results are now illustrated in the correctly specified case for the tvAR(1) model

(28) $$X_{t,n} + \alpha\left(\frac{t}{n}\right) X_{t-1,n} = \sigma\left(\frac{t}{n}\right) \epsilon_t$$

with independent Gaussian innovations $\epsilon_t$ satisfying $E\epsilon_t = 0$, $\text{Var}(\epsilon_t) = 1$, $\sigma(\cdot) > 0$ and $\sup_u |\alpha(u)| < 1$, with $\alpha(\cdot)$ smooth and $\sigma(\cdot)$ of bounded variation. These assumptions ensure that the corresponding time-varying spectral density $f(u, \lambda)$ exists and $\frac{1}{f(u,\lambda)} = \frac{2\pi}{\sigma^2(u)}(1 + \alpha^2(u) + 2\alpha(u)\cos(\lambda))$ [see (8)].

We will assume that $\alpha(\cdot) \in \mathcal{A}$ and $\sigma^2(\cdot) \in \mathcal{D}$, where $\mathcal{A}$ and $\mathcal{D}$ are model classes. This leads to

$$\mathcal{F}^* = AR(1; \mathcal{A}, \mathcal{D})$$
$$= \left\{ \frac{1}{g_{\alpha,\sigma^2}(u,\lambda)} = \frac{2\pi}{\sigma^2(u)}(1 + \alpha^2(u) + 2\alpha(u)\cos(\lambda)); \ \alpha \in \mathcal{A}; \ \sigma^2 \in \mathcal{D} \right\}.$$

PROOF. *Global estimation.* For simplicity, we assume here that $\sigma^2$ is a constant, that is, we choose

$$\mathcal{D} = \left( \epsilon^2, \frac{1}{\epsilon^2} \right)$$

for some $0 < \epsilon < 1$. We assume further that $\alpha$ is a member of the Sobolev space $\mathcal{S}^m$ with smoothness parameter $m \in \mathbf{N}$ such that the first $m \geq 1$ derivatives exist and have finite $L_2$-norms. To ensure $1 + \alpha(u)z \neq 0$ for all $0 < |z| \leq 1$ [cf. (6) ff.], we choose

$$\mathcal{A} = \{ h(\cdot) \in \mathcal{S}^m; \ \sup_u |h(u)| < 1 \}.$$

The metric entropy of $\mathcal{A}$ can be bounded by $A\eta^{-1/m}$ for some $A > 0$ [4]. It follows (under additional constraints on the model—see below for more details) that the metric entropy of the corresponding class of reciprocal spectral densities can be bounded by $\tilde{H}(\eta) = \tilde{A}\eta^{-\frac{1}{m}}$ for some $\tilde{A} > 0$. Hence, Theorem 2.6 gives us rates of convergence (by putting $\gamma = 1/m$). These rates are suboptimal and we can obtain faster rates via sieve estimation as we illustrate below.

*Sieve estimation.* Let $0 < \epsilon_n < 1$. We will assume $\epsilon_n \to 0$ as $n \to \infty$. We choose

$$\mathcal{D} = \mathcal{D}_n = \left( \epsilon_n^2, \frac{1}{\epsilon_n^2} \right).$$

For $0 < \epsilon_n < 1$ and $k_n$ a positive integer, let

$$\mathcal{A}_n = \left\{ \alpha(u) = a_0 + \sum_{j=1}^{k_n}(a_j \cos(2\pi j u) + b_j \sin(2\pi j u)); \right.$$
$$\left. u \in [0,1]; \sup_{u \in [0,1]} |\alpha(u)| < 1 \right\},$$

where $a_0, a_1, \ldots, a_{k_n}, b_1, \ldots, b_{k_n} \in \mathbf{R}$, and let

$$\mathcal{F}_n^* = AR(1; \mathcal{A}_n, \mathcal{D}_n).$$

It follows that $\sup_{u,\lambda} \phi(u,\lambda) = O(1/\epsilon_n^2)$ uniformly in $\phi \in \mathcal{F}_n^*$ and that $M^* = O(1/\epsilon_n^2)$. Note that we do not need the lower bound $M_*$: Kolmogorov's formula (cf. [5], Chapter 5.8) implies for all $u$ that $\int_{-\pi}^{\pi} \log g_{\alpha,\sigma^2}(u,\lambda)\,d\lambda = 2\pi \log(\sigma^2/(2\pi))$, leading to $R_{\log}(g_{\alpha,\sigma^2}) = 0$. As mentioned below Assumption 2.4, Theorem 2.6 can now be applied with $M_* = 1$. Further,

$$|\widehat{\phi}(u,0)| = \frac{2\pi}{\sigma^2}(1+\alpha^2(u)) \leq \frac{4\pi}{\epsilon_n^2},$$

$$|\widehat{\phi}(u,1)| = |\widehat{\phi}(u,-1)| = \frac{2\pi}{\sigma^2}|\alpha(u)| \leq \frac{2\pi}{\epsilon_n^2},$$

$$\widehat{\phi}(u,j) = 0 \quad \text{for } |j| \geq 2.$$

Consequently, $\rho_\infty$ and $v_\Sigma$ are of order $O(1/\epsilon_n^2)$ and it follows that $c_n = O((M^*)^2) = O(\frac{1}{\epsilon_n^4})$. As a finite-dimensional linear space of uniformly bounded functions, the metric entropy of $\mathcal{A}_n$ can, for small $\eta > 0$, be bounded by $A k_n \log(1/\epsilon_n \eta)$ and, hence, a similar upper bound of $\tilde{A} k_n \log(1/(\epsilon_n^2 \eta))$ holds for the metric entropy of $\mathcal{F}_n^*$. Finally, we determine the approximation error $a_n$. First, note that for $\epsilon_n \to 0$, we have $\sigma^2 \in \mathcal{D}_n$ for sufficiently large $n$. Further, for $\frac{1}{g_{\alpha_n,\sigma^2}} \in \mathcal{F}_n^*$ and $\frac{1}{g_{\alpha,\sigma^2}} \in \mathcal{F}^*$, we have $\rho_2(\frac{1}{g_{\alpha,\sigma^2}}, \frac{1}{g_{\alpha_n,\sigma^2}}) = O(\frac{1}{\sigma^2} \rho_2(\alpha, \alpha_n))$. It is well known that the approximation error of the sieve $\mathcal{A}_n$ in $\mathcal{S}^m$ is of the order $k_n^{-m}$ (e.g., see [3, 18, 21]). Hence, we can choose the approximating function $\pi_n(1/f_{\mathcal{F}})$ such that as $\epsilon_n \to 0$, we have

$$\rho_2\left(\frac{1}{f_{\mathcal{F}}}, \pi_n(1/f_{\mathcal{F}})\right) = O\left(\frac{1}{k_n^m}\right).$$

In other words, if $\epsilon_n \to 0$, we have $a_n = O(\frac{1}{k_n^m})$. We now choose the free parameters $k_n, \epsilon_n$ in order to balance the two terms in the definition of $\delta_n$. This leads us to the rate $\delta_n = \epsilon_n^{-4(1+\frac{m}{2m+1})}(\frac{n}{\log n})^{-\frac{m}{2m+1}}$. Choosing $\epsilon_n$ of the order $(\log n)^{-\alpha}$ for some $\alpha > 0$ gives us a rate which (up to a log term) equals the optimal rate known from the i.i.d. case. $\square$

Finally, we state that the same rates of convergence hold if the estimate is obtained by minimizing an approximation of the likelihood $\mathcal{L}_n(g)$. An example with a conditional likelihood function is given in Section 3.2.

THEOREM 2.8 (Likelihood approximations). *Let $\widetilde{\mathcal{L}}_n(g)$ be a criterion function with*

(29) $$\sup_{g \in \mathcal{F}_n} |\widetilde{\mathcal{L}}_n(g) - \mathcal{L}_n(g)| = o_P(\delta_n^2/(M^*)^2).$$

*Then Theorem 2.6 holds with $\widehat{f}_n$ replaced by $\widetilde{f}_n = \operatorname{argmin}_{g \in \mathcal{F}_n} \widetilde{\mathcal{L}}_n(g)$.*

2.4. *Discussion.* Why both sieve estimation and global estimation? The reason for considering both sieve and global estimation is more or less technical. In contrast to the standard empirical process, the (time-varying) empirical spectral process $E_n(\phi)$ is not monotonic in its argument $\phi$, that is, $\phi(u,\lambda) < \psi(u,\lambda) \,\forall u,\lambda$ does not imply $E_n(\phi) \leq E_n(\psi)$, since $J_n(u,\lambda)$ is not necessarily positive. This implies that the "bracketing idea" from standard empirical process theory cannot be applied. For this reason, we cannot fully exploit our Bernstein-type exponential inequality (36) below; essentially, we can only use the (less strong) exponential inequality (37). Therefore, we have to work with a covering integral which is the integral of the metric entropy instead of the square root of the metric entropy. As a consequence, we obtain slower rates of convergence. Our sieve estimators, however, do not suffer from this problem. At least if the model is correctly specified, then, as has been demonstrated in Example 2.7 and in Section 3 below, the sieve estimators achieve the same rates of convergence as the corresponding NPMLE of a probability density function based on i.i.d. data which, in this setting, are almost (i.e., up to log terms) optimal.

**3. Estimation under shape constraints.** Here, we consider the special case of a correctly specified model with constant AR-coefficients and monotonically increasing variance function. Our model spectral densities are hence of the form

$$g_{\alpha,\sigma^2}(\lambda) = \frac{\sigma^2(u)}{2\pi w_\alpha(\lambda)}, \qquad w_\alpha(\lambda) = \left|1 + \sum_{j=1}^p \alpha_j \exp(i\lambda j)\right|^2,$$

where the AR-coefficients $\alpha_1, \ldots, \alpha_p$ lie in the set

$$\mathcal{A}_p = \left\{\alpha = (\alpha_1, \ldots, \alpha_p)' \in \mathbf{R}^p : \left|1 + \sum_{j=1}^p \alpha_j z^j\right| \neq 0 \text{ for all } 0 < |z| \leq 1\right\},$$

[cf. (6)]. We assume that $\sigma^2(\cdot) \in \mathcal{M}$, where

$$\mathcal{M} = \left\{s^2 : [0,1] \to (0,\infty); \ s^2 \text{ increasing with } 0 < \inf_u s^2(u) \leq \sup_u s^2(u) < \infty\right\}.$$

With this notation, our model assumption can be formalized as

$$\mathcal{F}^* = AR(p; \mathcal{A}_p, \mathcal{M}) = \left\{\phi(u,\lambda) = \frac{2\pi}{\sigma^2(u)} w_\alpha(\lambda); \ \alpha \in \mathcal{A}_p, \ \sigma^2(\cdot) \in \mathcal{M}\right\}.$$

This model, with a unimodal rather than an increasing variance function, has been used in [6] for discriminating between earthquakes and explosions based on seismographic time series. To keep the exposition somewhat simpler, we shall only consider the special case of a monotonic instead of a unimodal variance.

*Global estimation.* Similarly to above, global estimation will lead to suboptimal rates of convergence. Since the class of bounded monotonic functions has a metric entropy satisfying $\log N(\eta, \mathcal{M}, \rho_2) \leq A\eta^{-1}$, and since the class $\mathcal{A}_p$ is finite-dimensional and hence its metric entropy is much smaller, it follows from Theorem 2.6 that our global NPMLE converges with rate $\delta_n = n^{-1/4}$ (provided all assumptions of this theorem are satisfied). As it turns out, this rate is suboptimal and can be improved upon by using sieve estimation. For this reason, we do not go into further detail concerning global estimation.

3.1. *Sieve estimation.* We first give a sieve for $\mathcal{M}$. For $k_n \in \mathbf{N}$ and $0 < \epsilon_n < 1$, let $\mathcal{C}_n = \mathcal{C}_n(\epsilon_n, k_n)$ denote the set of all increasing functions on $[0, 1]$, piecewise constant on the intervals $(\frac{j-1}{k_n}, \frac{j}{k_n}]$, $j = 1, \ldots, k_n$, and bounded from above and below by $1/\epsilon_n^2$ and $\epsilon_n^2$, respectively. Formally,

$$\mathcal{C}_n = \left\{ s^2(\cdot) : s^2(u) = \sum_{j=1}^{k_n} a_j \mathbf{1}\left(u \in \left(\frac{j-1}{k_n}, \frac{j}{k_n}\right]\right); \right.$$
$$\left. \epsilon_n^2 \leq a_1 \leq a_2 \leq \cdots \leq a_{k_n} \leq \frac{1}{\epsilon_n^2}, u \in [0,1] \right\}.$$

With these definitions, our sieve now becomes

$$\mathcal{F}_n^* = AR(p; \mathcal{A}_p, \mathcal{C}_n).$$

*Sieve estimation of the spectral density.* The next theorem states that we obtain with an appropriate choice of $\epsilon_n$ the known rate of $n^{-1/3}$ (up to a log term) for the NPMLE of the spectral density. This rate is known to be optimal for estimating a monotonic density based on i.i.d. data. Again, the proof is contained in Section 5.

THEOREM 3.1. *Let $X_{t,n}$ be a Gaussian locally stationary process and $\mathcal{F}^*$ and $\mathcal{F}_n^*$ be as defined above with $k_n = O(n^{1/3}(\log n)^{-2/3})$ and $\epsilon_n = (\log n)^{-1/5}$. If $f \in \mathcal{F}$, then we have*

$$\rho_2\left(\frac{1}{\widehat{f}_n}, \frac{1}{f}\right) = O_P(n^{-1/3} \log n).$$

*Sieve estimation of the monotonic variance function.* Next, we see that the above results for estimating the (inverse) spectral densities provide information about estimating the monotonic function $\sigma^2(\cdot)$ itself. We show that the rates of convergence from Theorem 3.1 translate to rates for $\widehat{\sigma}_n^2$. It can also be shown that the estimators of the finite-dimensional AR-parameters have a $\sqrt{n}$-rate and are asymptotically normal. This is not considered here, however.

Let
$$(\alpha_0, \sigma_0^2(\cdot)) = \operatorname*{argmin}_{(\alpha,\sigma) \in \mathcal{A}_p \times \mathcal{M}} \mathcal{L}(g_{\alpha,\sigma^2})$$

be the theoretically 'optimal' parameters and

$$(\widehat{\alpha}_n, \widehat{\sigma}_n^2(\cdot)) = \operatorname*{argmin}_{(\alpha,\sigma) \in \mathcal{A}_p \times \mathcal{C}_n} \mathcal{L}_n(g_{\alpha,\sigma^2})$$

be the sieve estimate.

THEOREM 3.2. *Let $X_{t,n}$ be a Gaussian locally stationary process and let $\mathcal{F}^*$ and $\mathcal{F}_n^*$ be as defined above, with $k_n = O(n^{1/3}(\log n)^{-2/3})$ and $\epsilon_n = (\log n)^{-1/5}$. If $f \in \mathcal{F}$, then*

$$(30) \qquad \rho_2\left(\frac{1}{\widehat{\sigma}_n^2}, \frac{1}{\sigma_0^2}\right) = O_P(n^{-1/3}\log n).$$

3.2. *An estimation algorithm.* Here, we discuss how to calculate a close approximation to the above $(\widehat{\alpha}_n, \widehat{\sigma}_n^2)$. The approximation considered is

$$(\widetilde{\alpha}_n, \widetilde{\sigma}_n^2(\cdot)) = \operatorname*{argmin}_{(\alpha,\sigma^2) \in \mathcal{A}_p \times \mathcal{C}_n} \widetilde{\mathcal{L}}_n(\alpha, \sigma^2),$$

where

$$(31) \quad \widetilde{\mathcal{L}}_n(\alpha, \sigma^2) = \frac{1}{n}\sum_{t=p+1}^{n}\left\{\log \sigma^2\left(\frac{t}{n}\right) + \frac{1}{\sigma^2(\frac{t}{n})}\left[X_{t,n} + \sum_{j=1}^{p}\alpha_j X_{t-j,n}\right]^2\right\}$$

is the so-called conditional Gaussian likelihood. By using Theorem 2.8, we now conclude that the minimizer has the same rate of convergence.

PROPOSITION 3.3. *Let $X_{t,n}$ be a Gaussian locally stationary process and let $\mathcal{F}^*$ and $\mathcal{F}_n^*$ be as defined above, with $k_n = O(n^{1/3}(\log n)^{-2/3})$ and $\epsilon_n = (\log n)^{-1/5}$. If $f \in \mathcal{F}$, then*

$$(32) \sup_{(\alpha,\sigma^2) \in \mathcal{A}_p \times \mathcal{C}_n} \left|\frac{1}{2}\{\widetilde{\mathcal{L}}_n(\alpha, \sigma^2) - \log(2\pi)\} - \mathcal{L}_n(g_{\alpha,\sigma^2})\right| = o_P(\delta_n^2/(M^*)^2).$$

*Hence, all assertions of Theorem 3.1 and Theorem 3.2 also hold for $\widetilde{f}_n = g_{\widetilde{\alpha}_n, \widetilde{\sigma}_n^2}$ and $\widetilde{\sigma}_n^2$, respectively.*

We now present our algorithm for calculating $(\widetilde{\alpha}_n, \widetilde{\sigma}_n^2(\cdot))$. Although global estimation is suboptimal, we first discuss the algorithm in this case in order to concentrate on the main ideas. The same ideas apply to sieve estimates, as will be indicated below.

Observe that for each fixed $\sigma^2(\cdot)$, minimizing $\widetilde{\mathcal{L}}_n(\alpha, \sigma^2)$ over $\alpha \in \mathcal{A}_p$ is a weighted least square problem. On the other hand, for each given $\alpha$, the minimizer over $\sigma^2(\cdot) \in \mathcal{M}$ can also be found explicitly. In fact, for each fixed $\alpha$, the minimizer

$$\widetilde{\sigma}^2_{n,\alpha}(\cdot) = \underset{\sigma^2 \in \mathcal{M}}{\operatorname{argmin}} \widetilde{\mathcal{L}}_n(\alpha, \sigma^2)$$

is given by the generalized isotonic regression to the squared residuals $e_t^2(\alpha) = (X_{t,n} + \sum_{j=1}^p \alpha_j X_{t-j,n})^2$. Note that there are no residuals for $t \leq p$ and, hence, the estimator is only defined for $t \geq p+1$. It follows that for $t \geq p+1$, the estimator $\widetilde{\sigma}^2_{n,\alpha}(\frac{t}{n})$ can be calculated as the (right) derivative of the greatest convex minorant to the cumulative sum diagram given by $\{(0,0), (\frac{t-p}{n-p}, \frac{1}{n-p} \sum_{s=p+1}^t e_s^2(\alpha)), t = p+1, \ldots, n\}$, by using the pool-adjacent-violators algorithm (PAVA). This follows from the theory of isotonic regression (cf. [20]). For completeness, let us briefly mention the relevant theory. Consider the expression

$$\tag{33} \sum_{t=p+1}^n (\Phi(x_t) - \Phi(y_t) - \phi(y_t)(x_t - y_t)),$$

where $\Phi$ is a convex function with derivative $\phi$. The theory of isotonic regression now implies that the minimizer of (33) over $(y_{p+1}, \ldots, y_n) \in \mathcal{K} = \{(y_{p+1}, \ldots, y_n) : y_{p+1} \leq \cdots \leq y_n\}$ is given by the (right) slope of the greatest convex minorant to the cumulative sums of the $x_t$'s, and it can be calculated by means of the PAVA. With $\Phi(x) = -\log x$, $x_t = e_t^2(\alpha)$ and $y_t = \sigma^2(t/n)$, we obtain

$$\Phi(x_t) - \Phi(y_t) - \phi(y_t)(x_t - y_t) = -\log e_t^2(\alpha) + \log \sigma^2(t/n) + \frac{e_t^2(\alpha) - \sigma^2(t/n)}{\sigma^2(t/n)}.$$

Consequently, for fixed $\alpha$, minimizing (33) over $\mathcal{K}$ is equivalent to minimizing $\widetilde{\mathcal{L}}_n(\alpha, \sigma^2)$ over all monotone functions $\sigma^2(\cdot) \in \mathcal{M}$. The global minimizer is then found by minimizing the profile likelihood $\widetilde{\mathcal{L}}_n(\alpha, \widetilde{\sigma}^2_{n,\alpha})$ over $\alpha$. Note that this is a continuous function in $\alpha$. This can be seen by observing that the squared residuals depend continuously on $\alpha$. Hence, at each fixed point $u = t/n$, the (right) slope of the greatest convex minorant, that is, $\widetilde{\sigma}^2_{n,\alpha}(t/n)$, is also a continuous function in $\alpha$. Therefore, we can conclude that the minimizer exists, provided the minimization is extended over a compact set of $\alpha$'s (as in the sieve estimation case).

*The basic algorithm* for finding the global minimizers $(\alpha, \sigma^2) \in \mathcal{A}_p \times \mathcal{M}$ is now given by the following iteration which results in a sequence $(\widetilde{\alpha}_{(k)}, \widetilde{\sigma}^2_{(k)})$, $k = 1, 2 \ldots$, with decreasing values of $\widetilde{\mathcal{L}}_n$. Given a starting value $\widetilde{\sigma}^2_{(0)}$, the iteration for $k = 1, 2, \ldots$ is as follows:

(i) Given $\widetilde{\sigma}^2_{(k-1)}(\cdot)$, find $\widetilde{\alpha}_k$ by solving the corresponding weighted least square problem

$$\widetilde{\alpha}_{(k)} = \underset{\alpha}{\operatorname{argmin}} \frac{1}{n} \sum_{t=p+1}^{n} \frac{1}{\widetilde{\sigma}^2_{(k-1)}(\frac{t}{n})} \left[ X_{t,n} + \sum_{j=1}^{p} \alpha_j X_{t-j,n} \right]^2.$$

(ii) Find $\widetilde{\sigma}^2_{(k)}(\cdot)$ as the solution to the PAVA using the squared residuals $e_t^2(\widetilde{\alpha}_{(k)})$, as outlined above.

A reasonable starting value $\widetilde{\sigma}^2_{(0)}(\cdot)$ is the solution of the PAVA using the squared raw data. This algorithm is applied in [6] to the problem of discrimination of time series.

The corresponding minimizer of the conditional Gaussian likelihood over the sieve parameter space $\mathcal{A}_p \times \mathcal{C}_n$ can be found similarly. First, note that the above solution is a piecewise constant function with jump locations in the set of points $\{\frac{t-p}{n-p}, t = p+1, \ldots, n\}$. Our sieve, however, consists of piecewise constant functions with jump locations in the set $\{\frac{j}{k_n}, j = 1, \ldots, k_n\}$. In order to find the minimizer of the conditional likelihood over this sieve, the only change one has to make to the above algorithm is to apply the PAVA to the cumulative sum diagram based on $\{(0,0), (\frac{t(j)-p}{n-p}, \frac{1}{n-p} \sum_{s=p+1}^{t(j)} e_s^2(\alpha)), j = \lceil \frac{k_n(p+1)}{n} \rceil, \ldots, k_n\}$, where $t(j) = \lceil \frac{nj}{k_n} \rceil$.

Note that in the above, we ignored the imposed boundedness restrictions on $\sigma^2$. An ad hoc way to obtain those is to truncate the isotonic regression at $1/\epsilon_n^2$ and from below at $\epsilon_n^2$. Alternatively, a solution respecting the bounds can be found by using the bounds $1/\epsilon_n^2$ and $\epsilon_n^2$ as upper and lower bounds in the PAVA. This means not to allow for derivatives outside this range, and to start the algorithm at $(0,0)$ and to end it at $(1, \frac{1}{n-p} \sum_{s=p+1}^{n} e_s^2(\widetilde{\alpha}_n))$. This can be achieved by a simple modification of the greatest convex minorant close to its endpoints. This also only modifies the estimator in the tail and close to the mode, but it has the additional advantage of sharing the property of the isotonic regression that the integral of (the piecewise constant function) $\widetilde{\sigma}_n^2$ equals the average of the squared residuals.

**4. The time-varying empirical spectral process.** In this section, we present exponential inequalities, maximal inequalities and weak convergence results for the empirical spectral process defined in (13). Asymptotic normality of the finite-dimensional distributions of $E_n(\phi)$ has been proved in [14]. Furthermore, several applications of the empirical spectral process have been discussed there. Let

(34) $$\rho_{2,n}(\phi) := \left( \frac{1}{n} \sum_{t=1}^{n} \int_{-\pi}^{\pi} \phi\left(\frac{t}{n}, \lambda\right)^2 d\lambda \right)^{1/2}$$

and

$$\tilde{E}_n(\phi) := \sqrt{n}(F_n(\phi) - \mathrm{E}F_n(\phi)). \tag{35}$$

THEOREM 4.1. *Let $X_{t,n}$ be a Gaussian locally stationary process, and let $\phi:[0,1]\times[-\pi,\pi]\to \mathbf{R}$ with $\rho_\infty(\phi)<\infty$ and $\tilde{v}(\phi)<\infty$. Then we have, for all $\eta>0$,*

$$P(|\tilde{E}_n(\phi)| \geq \eta) \leq c_1 \exp\left(-c_2 \frac{\eta^2}{\rho_{2,n}(\phi)^2 + \frac{\eta\rho_\infty(\phi)}{\sqrt{n}}}\right) \tag{36}$$

*and*

$$P(|\tilde{E}_n(\phi)| \geq \eta) \leq c_1 \exp\left(-c_2 \frac{\eta}{\rho_{2,n}(\phi)}\right) \tag{37}$$

*with some constants $c_1, c_2 > 0$. Furthermore, we have, for some $c_3 > 0$,*

$$\sqrt{n}|\mathrm{E}F_n(\phi) - F(\phi)| \leq c_3 n^{-1/2}(\rho_\infty(\phi) + \tilde{v}(\phi)). \tag{38}$$

REMARK 4.2. (i) Since $\rho_{2,n}(\phi)^2 \leq \rho_2(\phi)^2 + \frac{1}{n}\rho_\infty(\phi)\tilde{v}(\phi)$, we can replace $\rho_{2,n}(\phi)^2$ in (36) by the latter expression and $\rho_{2,n}(\phi)$ in (37) by $\rho_2(\phi) + (\frac{1}{n}\rho_\infty(\phi)\tilde{v}(\phi))^{1/2}$.

(ii) Combining (36) and (38) leads to the following exponential inequality for the empirical spectral process [see (66)]:

$$P(|E_n(\phi)| \geq \eta) \leq c_1' \exp\left(-c_2' \frac{\eta^2}{\rho_2(\phi)^2 + \frac{\eta}{n^{1/2}}(\rho_\infty(\phi) + \tilde{v}(\phi)) + \frac{1}{n}\rho_\infty(\phi)\tilde{v}(\phi)}\right).$$

However, we prefer to use the above inequalities and to treat the bias separately.

(iii) The constants $c_1, c_2, c_3$ depend on the characteristics of the process $X_{t,n}$, but not on $n$.

The above exponential inequalities form the core of the proofs of the next two results which lead to asymptotic stochastic equicontinuity of the empirical spectral process. Analogously to standard empirical process theory, stochastic equicontinuity is crucial for proving tightness. For proving the rates of convergence of the NPMLE, we need more, namely rates of convergence for the modulus of continuity. These rates also follow from the results below.

In the formulations of the following theorems, we use the constant

$$L = \max(c,1)^2 \max(K_1, K_2, K_3, 1) > 0, \tag{39}$$

where $c$ is the constant from Lemma A.4, the constants $K_1$–$K_3$ are from Lemma A.3 and $c_3$ is from (38). All of these constants do not depend on $n$ or on the function classes $\mathcal{F}_n$. They only depend on the constant K of the underlying process $X_{t,n}$ given in Definition 2.1.

THEOREM 4.3 (Sequence of "finite-dimensional" index classes). *Suppose that $X_{t,n}$ is a Gaussian locally stationary process. Further, let $\Phi_n$ be a function class satisfying*

$$\sup_{\phi \in \Phi_n} \rho_2(\phi) \leq \tau_2 < \infty. \tag{40}$$

*Assume further that there exist constants $A > 0$ and $k_n \geq 1$ such that*

$$\log N(\epsilon, \Phi_n, \rho_2) \leq A k_n \log\left(\frac{n}{\epsilon}\right) \qquad \forall \epsilon > 0.$$

*Let $d \geq 1$. Suppose that $\eta > 0$ satisfies the conditions*

$$\eta \leq d \frac{2 n^{\frac{1}{2}} \tau_2^2}{\rho_\infty}, \tag{41}$$

$$\eta \geq \tilde{c} n^{-\frac{1}{2}} \log n, \tag{42}$$

$$\eta^2 \geq \frac{24 A d}{c_2} \tau_2^2 k_n \log^+\left(\frac{8 L n^2 \log n}{\eta}\right), \tag{43}$$

*where $\tilde{c} \geq 2^4 L \max(v_\Sigma, 1)$. Then there exists a set $B_n$ with $\lim_{n \to \infty} P(B_n) = 1$, such that the inequality*

$$P\left(\sup_{\phi \in \Phi_n} |\tilde{E}_n(\phi)| > \eta; B_n\right) \leq c_1 \exp\left\{-\frac{c_2}{24 d} \frac{\eta^2}{\tau_2^2}\right\} \tag{44}$$

*holds, where $c_1, c_2 > 0$ are the constants from (36).*

The next result allows for richer model classes. It is formulated for fixed $n$ and, therefore, the class $\Phi$ may again depend on $n$. Since we apply this result for global estimation with a fixed class $\mathcal{F}^*$, it is formulated as if $\Phi$ were fixed.

THEOREM 4.4 ("Infinite-dimensional" index class). *Let $X_{t,n}$ be a Gaussian locally stationary process. Let $\Phi$ satisfy Assumption 2.4(c) (with $\Phi$ replacing $\mathcal{F}^*$) and suppose that (40) holds with $\tau_2 > 0$. Further, let $c_1$, $c_2$ be the constants from (37). There exists a set $B_n$ with $\lim_{n \to \infty} P(B_n) = 1$ such that*

$$P\left(\sup_{\phi \in \Phi} |\tilde{E}_n(\phi)| > \eta, B_n\right) \leq 3 c_1 \exp\left\{-\frac{c_2}{8} \frac{\eta}{\tau_2}\right\} \tag{45}$$

for all $\eta > 0$ satisfying the following conditions. Let $\alpha = \tilde{H}^{-1}(\frac{c_2}{8}\frac{\eta}{\tau_2})$. Then

(46) $$\eta \geq 2^6 L \max(\rho_\infty, v_\Sigma, 1) n^{-1/2} \log n$$

and either $\frac{\eta}{16 L n \log n} > \alpha$ or the following hold:

(47) $$\eta \geq 2\frac{\log 2}{c_2}\tau_2,$$

(48) $$\eta \geq \frac{192}{c_2} \int_{\frac{\eta}{2^6 L n \log n}}^{\alpha} \tilde{H}_\Phi(u)\,du.$$

By applying the above results to the class of differences $\{\phi_1 - \phi_2;\ \phi_1, \phi_2 \in \Phi_{(n)}, \rho_2(\phi_1 - \phi_2) \leq \delta_n\}$, we obtain rates of convergence for the modulus of continuity of the time-varying empirical spectral process. This is utilized in the proof of Theorem 2.6. As a special case, we obtain the asymptotic equicontinuity of $\{E_n(\phi), \phi \in \Phi\}$. Together with the convergence of the finite-dimensional distributions [14], this leads to the following functional central limit theorem (for weak convergence in $\ell^\infty(\Phi)$ we refer to [23], Section 1.5):

THEOREM 4.5. *Suppose that $X_{t,n}$ is a Gaussian locally stationary process. Let $\Phi$ be such that Assumption* 2.4(c) *holds (for $\Phi$ replacing $\mathcal{F}$) and*

(49) $$\int_0^1 \tilde{H}(u, \Phi, \rho_2)\,du < \infty.$$

*Then the process $(E_n(\phi); \phi \in \Phi)$ converges weakly in $\ell^\infty(\Phi)$ to a tight mean-zero Gaussian process $(E(\phi); \phi \in \Phi)$ with*

$$\operatorname{cov}(E(\phi_j), E(\phi_k))$$
$$= 2\pi \int_0^1 \frac{h^4(u)}{\|h\|_2^4} \int_{-\pi}^{\pi} \phi_j(u,\lambda)[\phi_k(u,\lambda) + \phi_k(u,-\lambda)]f^2(u,\lambda)\,d\lambda\,du.$$

**5. Proofs.** The following lemma establishes the relation between the $L_2$-norm and the Kullback–Leibler information divergence. This relation is a key ingredient in the proof of Theorem 2.6.

LEMMA 5.1. *Let $\mathcal{F}$ be such that $f_\mathcal{F}$ exists and is unique and $\mathcal{L}(g) < \infty$ for all $g \in \mathcal{F}$.*

(i) *Assume that for some constant $0 < M^* < \infty$, we have $\sup_{u,\lambda} |\phi(u,\lambda)| < M^*$ for all $\phi \in \mathcal{F}^*$. If $\mathcal{F}^*$ is convex, then for all $g \in \mathcal{F}$,*

$$\frac{1}{8\pi(M^*)^2}\rho_2^2\left(\frac{1}{g}, \frac{1}{f_\mathcal{F}}\right) \leq \mathcal{L}(g) - \mathcal{L}(f_\mathcal{F}).$$

*If the model is correctly specified, that is, $f_\mathcal{F} = f$, then the above inequality holds without the convexity assumption.*

(ii) *Assume that for some constant $0 < M_* < \infty$, we have $M_* < \inf_{u,\lambda} |\phi(u,\lambda)|$ for all $\phi \in \mathcal{F}^*$. Then we have, with $\max(\sup_{u,\lambda} |f(u,\lambda)|, 1/M_*) < \Omega < \infty$ for all $g \in \mathcal{F}$,*

$$\mathcal{L}(g) - \mathcal{L}(f) \leq \frac{\Omega^2}{4\pi} \rho_2^2\left(\frac{1}{g}, \frac{1}{f}\right) \quad \text{and} \tag{50}$$

$$\mathcal{L}(g) - \mathcal{L}(f_\mathcal{F}) \leq \frac{\Omega}{2\pi} \max\left(\Omega \rho_2^2\left(\frac{1}{g}, \frac{1}{f_\mathcal{F}}\right), \rho_1\left(\frac{1}{g}, \frac{1}{f_\mathcal{F}}\right)\right), \tag{51}$$

*where $\rho_1$ denotes the $L_1$-norm on $[0,1] \times [-\pi, \pi]$.*

PROOF. For $g \in \mathcal{F}$, let $w = 1/g \in \mathcal{F}^*$. We set

$$\Psi(w) = \mathcal{L}\left(\frac{1}{w}\right) = \frac{1}{4\pi} \int_0^1 \int_{-\pi}^{\pi} \{-\log w(u,\lambda) + w(u,\lambda) f(u,\lambda)\} \, d\lambda \, du.$$

Direct calculations yield the following Gateaux derivative $\delta\Psi$ of $\Psi$: For $v, w \in \mathcal{F}^*$,

$$\delta\Psi(v, w) := \frac{\partial}{\partial t} \Psi(v + tw)|_{t=0} = \frac{1}{4\pi} \int_0^1 \int_{-\pi}^{\pi} \left\{-\frac{w}{v} + wf\right\} d\lambda \, du. \tag{52}$$

It follows that

$$\Psi(w) - \Psi(v) - \delta\Psi(v, w - v) = \frac{1}{4\pi} \int_0^1 \int_{-\pi}^{\pi} \left\{-\log \frac{w}{v} + \frac{w}{v} - 1\right\} d\lambda \, du. \tag{53}$$

Since $v = v(u,\lambda)$ and $w = w(u,\lambda)$ are uniformly bounded and $\log(1+x) = x + R(x)$, where $R(x) = -\frac{x^2}{2(1+\theta x)^2}$ with some $\theta \in (0,1)$, we obtain uniformly in $u$ and $\lambda$,

$$-\log \frac{w}{v} + \frac{w}{v} - 1 = \frac{1}{2}(w-v)^2 / (v + \theta(w-v))^2 \geq \frac{1}{2(M^*)^2}(w-v)^2. \tag{54}$$

Hence,

$$\Psi(w) - \Psi(v) - \delta\Psi(v, w-v) \geq \frac{1}{8\pi(M^*)^2} \rho_2(w,v)^2. \tag{55}$$

Therefore, the function $\Psi$ is strongly convex on $\mathcal{F}^*$. Corollary 10.3.3 of Eggermont and LaRiccia [15] now implies the inequality from (i), provided $\mathcal{F}^*$ is convex. Note further that if the model is correct, that is, if $f_\mathcal{F} = f$, then it is straightforward to see that (53) holds for $v = 1/f$ if we formally let $\delta\Psi(\frac{1}{f}, w - \frac{1}{f}) = 0$. In other words, if the model is correctly specified, the result follows directly from (55) without using the convexity assumption.

As for the second part of the lemma, observe that similarly to (54), the assumed boundedness of $g$ and $f_\mathcal{F}$ implies that

$$\frac{1}{4\pi} \int_0^1 \int_{-\pi}^{\pi} \left\{-\log \frac{f_\mathcal{F}}{g} + \frac{f_\mathcal{F}}{g} - 1\right\} du \, d\lambda \leq \frac{\Omega^2}{8\pi} \rho^2\left(\frac{1}{g}, \frac{1}{f_\mathcal{F}}\right)^2.$$

Further, $\delta\Psi(\frac{1}{f_\mathcal{F}}, w - \frac{1}{f_\mathcal{F}}) = \frac{1}{4\pi}\int_0^1\int_{-\pi}^{\pi}(f - f_\mathcal{F})(w - 1/f_\mathcal{F})\,du\,d\lambda$. Hence,

$$\left|\delta\Psi\left(\frac{1}{f_\mathcal{F}}, w - \frac{1}{f_\mathcal{F}}\right)\right| \leq \frac{\Omega}{4\pi}\rho_1(w, 1/f_\mathcal{F}).$$

Notice that if $f = f_\mathcal{F}$, then $\delta\Psi(\frac{1}{f_\mathcal{F}}, w - \frac{1}{f_\mathcal{F}}) = 0$. The result follows by using (53). $\square$

PROOF OF THEOREM 2.6. First we prove part II. Lemma 5.1 and (18), along with the fact $\tilde{E}_n(\phi) = \sqrt{n}(F_n(\phi) - \mathrm{E}F_n(\phi))$ lead to the relation

$$
\begin{aligned}
P\bigg[\rho_2\bigg(\frac{1}{\widehat{f}_n}, \frac{1}{f_\mathcal{F}}\bigg) &\geq C\delta_n\bigg] \\
&\leq P\bigg[\sup_{g\in\mathcal{F};\rho_2(\frac{1}{g},\frac{1}{f_\mathcal{F}})\geq C\delta_n}\frac{\tilde{E}_n(\frac{1}{f_\mathcal{F}} - \frac{1}{g})}{\rho_2^2(\frac{1}{g},\frac{1}{f_\mathcal{F}})} \geq \frac{\sqrt{n}}{4(M^*)^2}\bigg] \\
&\quad + P\bigg[R_n \geq \frac{C^2\delta_n^2}{16\pi(M^*)^2}\bigg],
\end{aligned}
\tag{56}
$$

where

$$R_n = \frac{1}{4\pi}(\mathrm{E}F_n - F)\bigg(\frac{1}{f_\mathcal{F}} - \frac{1}{\widehat{f}_n}\bigg) + R_{\log}(f_\mathcal{F}) - R_{\log}(\widehat{f}_n).$$

Note that the expectation operator E only operates on $F_n$ and not on $\widehat{f}_n$. Theorem 4.1, Lemma A.2 and (38) imply that the second term in (56) tends to zero as $n \to \infty$.

We now use the so-called "peeling device" (e.g., [22]) in combination with Theorem 4.4. The first term in (56) is bounded by

$$
\begin{aligned}
P\bigg[&\sup_{2\rho_\infty \geq \rho_2(\phi,\psi) \geq C\delta_n}\frac{|\tilde{E}_n(\phi - \psi)|}{\rho_2^2(\phi,\psi)} \geq \frac{\sqrt{n}}{4(M^*)^2}\bigg] \\
&\leq \sum_{j=0}^{K_n} P\bigg[\sup_{C2^{j+1}\delta_n \geq \rho_2(\phi,\psi) \geq C2^j\delta_n}\frac{|\tilde{E}_n(\phi - \psi)|}{\rho_2^2(\phi,\psi)} \geq \frac{\sqrt{n}}{4(M^*)^2}\bigg] \\
&\leq \sum_{j=0}^{K_n} P\bigg[\sup_{\rho_2(\phi,\psi) \leq C2^{j+1}\delta_n}|\tilde{E}_n(\phi - \psi)| \geq \frac{\sqrt{n}(C2^j\delta_n)^2}{4(M^*)^2}\bigg],
\end{aligned}
\tag{57}
$$

where $K_n$ is such that $\lfloor C2^{K_n+1}\delta_n\rfloor = 2\rho_\infty$. For sufficiently large $C$, our assumptions ensure that the conditions of Theorem 4.4 are satisfied for all $j = 1, \ldots, K_n$, [with $\eta = \sqrt{n}(C2^j\delta_n)^2/4(M^*)^2$ and $\tau_2 = C2^{j+1}\delta_n$ when $\delta_n$ is chosen as stated in the theorem we are proving]. Hence, on an exceptional set whose probability tends to zero as $n \to \infty$, we obtain the bound

$$\leq 3\sum_{j=0}^{\infty}c_1\exp\bigg\{-c_2\frac{\sqrt{n}C2^j\delta_n}{64(M^*)^2}\bigg\},$$

which tends to zero as $C \to \infty$ because $n\delta_n^2$, by assumption, is bounded away from zero.

The proof of part $I$ has a similar structure. First, we derive the analog to (18). Let $\pi_n = \pi_n(1/f_{\mathcal{F}}) \in \mathcal{F}_n^*$ denote an approximation to $1/f_{\mathcal{F}}$ in $\mathcal{F}_n^*$ (cf. Assumption 2.4). Then

$$
\begin{aligned}
(58) \quad 0 &\leq \mathcal{L}_n(1/\pi_n) - \mathcal{L}_n(\widehat{f}_n) \\
&= (\mathcal{L}_n(1/\pi_n) - \mathcal{L}(1/\pi_n)) - (\mathcal{L}_n(\widehat{f}_n) - \mathcal{L}(\widehat{f}_n)) - (\mathcal{L}(\widehat{f}_n) - \mathcal{L}(f_{\mathcal{F}})) \\
&\quad - (\mathcal{L}(f_{\mathcal{F}}) - \mathcal{L}(1/\pi_n)).
\end{aligned}
$$

In the following, we prove the case of a correctly specified model, that is, $f = f_{\mathcal{F}}$. Using (50) from Lemma 5.1, we obtain

$$
\begin{aligned}
&\frac{1}{8\pi(M^*)^2}\rho_2^2\left(\frac{1}{\widehat{f}_n}, \frac{1}{f_{\mathcal{F}}}\right) \\
&\quad \leq (\mathcal{L}_n(1/\pi_n) - \mathcal{L}(1/\pi_n)) - (\mathcal{L}_n(\widehat{f}_n) - \mathcal{L}(\widehat{f}_n)) + (\mathcal{L}(1/\pi_n) - \mathcal{L}(f_{\mathcal{F}})) \\
(59) \quad &\quad = \frac{1}{4\pi}\frac{1}{\sqrt{n}}E_n\left(\pi_n - \frac{1}{\widehat{f}_n}\right) + R_{\log}(1/\pi_n) - R_{\log}(\widehat{f}_n) + (\mathcal{L}(1/\pi_n) - \mathcal{L}(f_{\mathcal{F}})) \\
&\quad \leq \frac{1}{4\pi}\sup_{g \in \mathcal{F}_n}(F_n - F)\left(\pi_n - \frac{1}{g}\right) + 2\sup_{g \in \mathcal{F}_n}|R_{\log}(g)| + \frac{\Omega^2}{2\pi}\rho_2^2(\pi_n, 1/f_{\mathcal{F}}).
\end{aligned}
$$

Note that $\Omega = \max(1/M_*, m)$, where $m = \sup_{u,\lambda} f(u,\lambda)$. Hence, $\Omega = 1/M_* \times \max(1, M_* m) \leq c/M_*$ with $c = \max(1, m)$. Let $\pi_n$ be a "good" approximation to $1/f_{\mathcal{F}}$ in $\mathcal{F}^*$, in the sense that $\rho_2^2(\pi_n, 1/f_{\mathcal{F}}) \leq 2a_n^2$. It follows that on the set $\{\rho_2^2(\frac{1}{\widehat{f}_n}, \frac{1}{f_{\mathcal{F}}}) \geq (C^2+2)\delta_n^2\}$, we have, by definition of $\delta_n$ and by using the triangle inequality, that $\rho_2^2(\frac{1}{\widehat{f}_n}, \pi_n) \geq \rho_2^2(\frac{1}{\widehat{f}_n}, \frac{1}{f_{\mathcal{F}}}) - \rho_2^2(\pi_n, \frac{1}{f_{\mathcal{F}}}) \geq (C^2+2)\delta_n^2 - 2a_n^2 \geq C^2\delta_n^2$. Hence, we obtain, by utilizing (59),

$$
\begin{aligned}
&P\left[\rho_2^2\left(\frac{1}{\widehat{f}_n}, \frac{1}{f_{\mathcal{F}}}\right) \geq (C^2+2)\delta_n^2\right] \\
&\quad \leq P\left[\sup_{g \in \mathcal{F}_n; \rho_2^2(\frac{1}{g}, \pi_n) \geq C\delta_n} \frac{\tilde{E}_n(\pi_n - \frac{1}{g})}{\rho_2^2(\frac{1}{g}, \pi_n)} \geq \frac{\sqrt{n}}{4(M^*)^2}\right] + P\left[R_n \geq \frac{(C^2+2)\delta_n^2}{32\pi(M^*)^2}\right] \\
&\quad\quad + P\left[\rho_2^2(\pi_n, 1/f_{\mathcal{F}}) > \frac{M_*^2(C^2+2)\delta_n^2}{16c(M^*)^2}\right].
\end{aligned}
$$

By definition of $\delta_n$, the last term is zero for $C$ sufficiently large. At this point, the proof is completely analogous to the proof of part $II$. The same arguments as used above show that the second term on the right-hand side can be made arbitrarily small by choosing $C$ sufficiently large. To bound the first term in the last inequality, we use the peeling device as above and Theorem 4.3 [with $d = \max(\frac{\rho_\infty}{(M^*)^2}, 1)$] to bound the sum analogously to (57).

We end up with the bound

$$\leq \sum_{j=0}^{\infty} c_1 \exp\biggl\{-\frac{c(C 2^j \delta_n)^2 n}{\rho_\infty(M^*)^2}\biggr\},$$

for some constant $c > 0$ and this sum also tends to 0 as $C \to \infty$.

The proof for $f \neq f_\mathcal{F}$ is, mutatis mutandis, the same. Instead of (58), we start with a version of (18) where $1/f$ is replaced by $1/f_{\mathcal{F}_n}$. We then proceed analogously to the proof above. □

PROOF OF THEOREM 2.8. As in (18), we obtain

$$0 \leq \mathcal{L}(\widetilde{f}_n) - \mathcal{L}(f_\mathcal{F}) \leq (\widetilde{\mathcal{L}}_n - \mathcal{L})(f_\mathcal{F}) - (\widetilde{\mathcal{L}}_n - \mathcal{L})(\widetilde{f}_n)$$

$$\leq (\mathcal{L}_n - \mathcal{L})(f_\mathcal{F}) - (\mathcal{L}_n - \mathcal{L})(\widetilde{f}_n) + o_P(\delta_n^2/(M^*)^2)$$

and, therefore, the proof of Theorem 2.6 applies, where $R_n$ is replaced by

$$\widetilde{R}_n = \frac{1}{4\pi}(\mathrm{E}F_n - F)\biggl(\frac{1}{f_\mathcal{F}} - \frac{1}{\widetilde{f}_n}\biggr) + R_{\log}(f_\mathcal{F}) - R_{\log}(\widetilde{f}_n) + o_P(\delta_n^2/(M^*)^2). \quad \square$$

PROOF OF THEOREM 3.1. The proof is an application of Theorem 2.6. We first derive the approximation error. For an increasing function $\sigma^2 \in \mathcal{M}$ with $\epsilon \leq \sigma^2(\cdot) \leq b$, let $s_n^2 \in \mathcal{C}_n$ with $a_j = s^2(j/k_n)$. Clearly, if $1/\epsilon_n^2 > b$ and $\epsilon_n^2 < \epsilon$, then

$$\sqrt{\int_0^1 (\sigma^2(u) - s_n^2(u))^2 \, du} \leq b/k_n.$$

In other words, the approximation error of $\mathcal{C}_n$ as a sieve for $\mathcal{M}$ is $O(1/k_n)$, provided $\epsilon_n \to 0$, which implies that $a_n = O(\frac{1}{k_n})$. Next, we determine a bound on the metric entropy of $\mathcal{F}_n^*$. Observe that as a space of functions bounded by $1/\epsilon_n^2$ and spanned by $k_n$ functions, the metric entropy of $\mathcal{C}_n$ satisfies $\log N(\eta, \mathcal{C}_n, \rho_2) \leq A k_n \log(1/(\epsilon_n^2 \eta))$ for some $A > 0$ (e.g., [22], Corollary 2.6). Next, we derive a bound for the metric entropy of $W_p = \{w_\alpha; \alpha \in \mathcal{A}_p\}$. First, note that $\log N(\eta, \mathcal{A}_p, \rho_2) \leq A \log(1/\eta)$ for some $A > 0$ (since $\mathcal{A}_p$ is bounded). Since $w_\alpha(\lambda) \leq 2^{2p}$ and

$$|w_{\alpha_1}(\lambda) - w_{\alpha_2}(\lambda)| \leq 2^p \sum_{j=1}^p |\alpha_{1j} - \alpha_{2j}|,$$

this leads to the bound $\log N(\eta, W_p, \rho_2) \leq \tilde{A} \log(1/\eta)$ for some $\tilde{A} > 0$. The two bounds on the metric entropy of $\mathcal{C}_n$ and $W_p$ now translate into the bound $\log N(\eta, \mathcal{F}_n^*, \rho_2) \leq \tilde{A} k_n \log(1/(\epsilon_n^2 \eta))$ for some $\tilde{A} > 0$. This can be seen as follows. First notice that with $s_n^2$, the approximation on $\sigma^2$ in $\mathcal{C}_n$ defined above, we have

$$\rho_2\left(\frac{1}{\sigma^2}w_{\alpha_1} - \frac{1}{s_n^2}w_{\alpha_2}\right) \leq \rho_2\left(\left(\frac{1}{\sigma^2} - \frac{1}{s_n^2}\right)w_{\alpha_1}\right) + \rho_2\left((w_{\alpha_1} - w_{\alpha_2})\frac{1}{s_n^2}\right).$$

Observe further that $\frac{1}{\sigma^2} \in \mathcal{C}_n^{-1} := \{\frac{1}{s^2}; s^2 \in \mathcal{C}_n\}$ and note that as a class of functions, $\mathcal{C}_n^{-1}$ is very similar to the class $\mathcal{C}_n$. The only difference is that $\mathcal{C}_n^{-1}$ consists of decreasing functions rather than increasing functions. In particular, the bound given above for the metric entropy of $\mathcal{C}_n$ also applies to the metric entropy of $\mathcal{C}_n^{-1}$. Since, in addition, $\frac{1}{s_n^2} \leq \frac{1}{\epsilon_n^2}$ and $w_\alpha \leq 1$, one sees that

$$N(\eta, \mathcal{F}_n^*, \rho_2) \leq N(\eta, \mathcal{C}_n^{-1}, \rho_2) N(\eta \epsilon_n^2, W_p, \rho_2).$$

This leads to the asserted bound for the metric entropy.

Note further that we can choose $M^* = O(1/\epsilon_n^2)$. As in Example 2.7, we avoid the lower bound $M_*$ by looking at $R_{\log}(g_{\alpha,\sigma^2})$ separately: We have for all $u$ that $\int_{-\pi}^{\pi} \log g_{\alpha,\sigma^2}(u,\lambda)\,d\lambda = 2\pi \log(\sigma^2(u)/(2\pi))$ and, therefore,

$$\sup_{g \in \mathcal{F}_n} |R_{\log}(g)| = O\left(\frac{\log(1/\epsilon_n)}{n}\right),$$

which is sufficient in our situation [cf. (56)]. Hence, the 'best' rate for the NPMLE which can be obtained from (26) follows by balancing $\sqrt{\frac{c_n k_n \log n}{n}}$ and $a_n$, where here, $c_n = O(\frac{1}{\epsilon_n^4})$. The latter follows as in Example 2.7. In the correctly specified case (i.e., the variance function actually is monotonic), this gives $k_n = (\frac{n\epsilon_n^4}{\log n})^{1/3}$ and, hence, the rate $\delta_n = \epsilon_n^{-10/3}(\frac{n}{\log n})^{-1/3}$. If we choose $1/\epsilon_n = O((\log n)^{1/5})$, then the rate becomes $n^{-1/3}\log n$, as asserted. □

PROOF OF THEOREM 3.2. For the true spectral density $f(u,\lambda) = s^2(u)/(2\pi v(\lambda))$, we assume, without loss of generality, that $\int_{-\pi}^{\pi} \log v(\lambda)\,d\lambda = 0$. This can be achieved by multiplying $v(\cdot)$ by an adequate constant. Since $\int_{-\pi}^{\pi} \log w_\alpha(\lambda)\,d\lambda = 0$ (Kolmogorov's formula) and $(\alpha_0, \sigma_0^2(\cdot))$ minimizes $\mathcal{L}(g_{\alpha,\sigma^2})$ over $\mathcal{A}_p \times \mathcal{M}$, we have

$$0 \leq \mathcal{L}(g_{\widehat{\alpha},\sigma_0^2}) - \mathcal{L}(g_{\alpha_0,\sigma_0^2}) = \frac{1}{4\pi}\int_0^1 \frac{s^2(u)}{\sigma_0^2(u)}\,du \int_{-\pi}^{\pi}\{w_{\widehat{\alpha}}(\lambda) - w_{\alpha_0}(\lambda)\}\frac{1}{v(\lambda)}\,d\lambda,$$

that is,

$$\int_{-\pi}^{\pi} \frac{w_{\widehat{\alpha}}(\lambda)}{v(\lambda)}\,d\lambda \geq \int_{-\pi}^{\pi} \frac{w_{\alpha_0}(\lambda)}{v(\lambda)}\,d\lambda.$$

Hence, we have

$$\mathcal{L}(g_{\widehat{\alpha},\widehat{\sigma}^2}) - \mathcal{L}(g_{\alpha_0,\sigma_0^2})$$
$$= \frac{1}{4\pi} \int_0^1 \left[ \log \frac{\widehat{\sigma}^2(u)}{\sigma_0^2(u)} + \int_{-\pi}^{\pi} \left\{ \frac{s^2(u)w_{\widehat{\alpha}}(\lambda)}{\widehat{\sigma}^2(u)v(\lambda)} - \frac{s^2(u)w_{\alpha_0}(\lambda)}{\sigma_0^2(u)v(\lambda)} \right\} d\lambda \right] du$$
$$\geq \frac{1}{4\pi} \int_0^1 \left[ \log \frac{\widehat{\sigma}^2(u)}{\sigma_0^2(u)} + \int_{-\pi}^{\pi} \left\{ \frac{s^2(u)}{\widehat{\sigma}^2(u)} - \frac{s^2(u)}{\sigma_0^2(u)} \right\} \frac{w_{\alpha_0}(\lambda)}{v(\lambda)} d\lambda \right] du$$
$$= \mathcal{L}(g_{\alpha_0,\widehat{\sigma}^2}) - \mathcal{L}(g_{\alpha_0,\sigma_0^2}).$$

For $s^2 \in \mathcal{C}_n$, let $\mathcal{H}(s^2) = \mathcal{L}(g_{\alpha_0,s^2})$. Then, obviously, $\sigma_0^2 = \operatorname{argmin}_{s^2 \in \mathcal{C}_n} \mathcal{H}(s^2)$. As in the proof of Lemma 5.1, it follows that

(60) $\quad \mathcal{L}(g_{\alpha_0,\widehat{\sigma}^2}) - \mathcal{L}(g_{\alpha_0,\sigma_0^2}) = \mathcal{H}(\widehat{\sigma}^2) - \mathcal{H}(\sigma_0^2) \geq \frac{1}{8\pi} \epsilon_n^4 \rho_2 \left( \frac{1}{\widehat{\sigma}^2}, \frac{1}{\sigma_0^2} \right)^2.$

Here, we use the fact that an upper bound for the functions in $\mathcal{C}_n$ is given by $1/\epsilon_n^2$. Assertion (30) now follows since we know from Theorem 3.1, the proof of Theorem 2.6 and Lemma 5.1 that

$$\frac{1}{\epsilon_n^4}(\mathcal{L}(g_{\widehat{\alpha},\widehat{\sigma}^2}) - \mathcal{L}(g_{\alpha_0,\sigma_0^2})) = O_P\left(\rho_2^2\left(\frac{1}{\widehat{f}_n}, \frac{1}{f_{\mathcal{F}}}\right)\right) = O_P(n^{-2/3}(\log n)^2).$$

Here, we used the fact that an upper bound for the functions in $\mathcal{F}^*$ is given by $\frac{1}{\epsilon_n^2}$ and the fact that with our choice of $\epsilon_n$, the rate of convergence for the NPMLE for the spectral density is $O_P(n^{-1/3} \log n)$ (cf. Theorem 3.1). $\square$

PROOF OF PROPOSITION 3.3. We obtain, with (8)–(10) and Kolmogorov's formula,

$$\mathcal{L}_n(g_{\alpha,\sigma^2}) + \frac{1}{2}\log(2\pi)$$
$$= \frac{1}{2n} \sum_{t=1}^n \log \sigma^2\left(\frac{t}{n}\right)$$
$$+ \frac{1}{2n} \sum_{t=1}^n \frac{1}{\sigma^2(\frac{t}{n})} \sum_{j,k=0}^p \alpha_j \alpha_k X_{[t+1/2+(j-k)/2],n}$$
$$\times X_{[t+1/2-(j-k)/2],n} \mathbf{1}_{1 \leq [t+1/2\pm(j-k)/2] \leq n},$$

where $\alpha_0 = 1$. The second summand is equal to

$$\frac{1}{2n} \sum_{j,k=0}^p \alpha_j \alpha_k \sum_{\{t:\, 1 \leq t-j, t-k \leq n\}} \frac{1}{\sigma^2\left(\frac{[t-j/2-k/2]}{n}\right)} X_{t-j,n} X_{t-k,n}.$$

By using the definition of $\widetilde{\mathcal{L}}_n(\alpha, \sigma^2)$ in (31), Lemma A.4 and the monotonicity of $\sigma^2(\cdot)$, we therefore obtain, with $\delta_n = n^{-1/3} \log n$ and $M^* = 1/\epsilon_n^2$,

$$\sup_{(\alpha,\sigma^2) \in \mathcal{A}_p \times \mathcal{C}_n} \left| \frac{1}{2} \{\widetilde{\mathcal{L}}_n(\alpha, \sigma^2) - \log(2\pi)\} - \mathcal{L}_n(g_{\alpha,\sigma^2}) \right| = O_p\left(\frac{\log(1/\epsilon_n)}{n} + \frac{\log n}{n\epsilon_n^2}\right)$$
$$= o_P(\delta_n^2/(M^*)^2).$$

PROOF OF THEOREM 4.1. We start with two technical lemmas. Direct calculation shows that

(61) $$F_n(\phi) = \frac{1}{n} \underline{X}'_n U_n\left(\frac{1}{2\pi}\phi\right) \underline{X}_n,$$

where $U_n(\phi)_{jk} = \hat{\phi}(\frac{1}{n}\lfloor \frac{j+k}{2} \rfloor, j - k)$ and $\lfloor x \rfloor$ denotes the largest integer less than or equal to $x$. The properties of $U_n(\phi)$ have been investigated under different assumptions in [10]. In this paper, we only need the following result on $\|U_n(\phi)\|_{\mathrm{spec}}$ and $\|U_n(\phi)\|_2$, where $\|A\|_2 := \mathrm{tr}(A^H A)^{1/2} = (\sum_{i,j} |a_{ij}|^2)^{1/2}$ is the Euclidean norm and $\|A\|_{\mathrm{spec}} := \sup_{\|x\|_2=1} \|Ax\|_2 = \max\{\sqrt{\lambda} \mid \lambda \text{ eigenvalue of } A^H A\}$ is the spectral norm:

LEMMA 5.2. *With $\rho_2(\phi)$, $\rho_{2,n}(\phi)$, $\rho_\infty(\phi)$ and $\tilde{v}(\phi)$ as defined in* (21), (34) *and* (23), *we have*

(62) $$\|U_n(\phi)\|_{\mathrm{spec}} \leq \rho_\infty(\phi)$$

*and*

(63) $$n^{-1}\|U_n(\phi)\|_2^2 \leq 2\pi \rho_{2,n}(\phi)^2 \leq 2\pi \rho_2(\phi)^2 + \frac{2\pi}{n}\rho_\infty(\phi)\tilde{v}(\phi).$$

PROOF. Let

$$\hat{\phi}_{jk} := \hat{\phi}\left(\frac{1}{n}\left\lfloor \frac{j+k}{2} \right\rfloor, j - k\right).$$

Then for $x \in \mathbf{C}^n$ with $\|x\|_2 = 1$,

(64) $$\|U_n(\phi)x\|_2^2 = \sum_{i,j,k=1}^n \bar{x}_i \hat{\phi}_{ji} \hat{\phi}_{jk} x_k = \sum_{j,\ell,m} \bar{x}_{j+\ell} \hat{\phi}_{j,j+\ell} \hat{\phi}_{j,j+m} x_{j+m}$$
$$\leq \sum_{\ell,m} \sup_j |\hat{\phi}_{j,j+\ell}| \sup_j |\hat{\phi}_{j,j+m}| \sum_j |\bar{x}_{j+\ell} x_{j+m}|,$$

where the range of summation is such that $1 \leq j + \ell$, $j + m \leq n$. Since $\Sigma_j |\bar{x}_{j+\ell} x_{j+m}| \leq \|x\|_2^2 = 1$, we obtain the first part. Furthermore, we have, with Parseval's equality,

$$\frac{1}{n}\|U_n(\phi)\|_2^2 = \frac{1}{n}\sum_{j,k=1}^{n}|\hat{\phi}_{jk}|^2 \leq \frac{1}{n}\sum_{t=1}^{n}\sum_{\ell=-\infty}^{\infty}\left|\hat{\phi}\left(\frac{t}{n},\ell\right)\right|^2 = 2\pi\rho_{2,n}(\phi)^2$$

(65)
$$= \int_0^1 \sum_{\ell=-\infty}^{\infty}|\hat{\phi}(u,\ell)|^2 du$$
$$+ \sum_{t=1}^{n}\int_0^{1/n}\sum_{\ell=-\infty}^{\infty}\left\{\left|\hat{\phi}\left(\frac{t}{n},\ell\right)\right|^2 - \left|\hat{\phi}\left(\frac{t-1}{n}+x,\ell\right)\right|^2\right\}dx$$
$$\leq 2\pi\rho_2(\phi)^2 + \frac{2\pi}{n}\rho_\infty(\phi)\tilde{v}(\phi). \qquad \square$$

LEMMA 5.3.  *If $\Sigma_n$ is the covariance matrix of the random vector $(X_{1,n},\ldots,X_{n,n})'$, then*

$$\|\Sigma_n^{1/2}\|_{\mathrm{spec}}^2 \leq \sum_{k=-(n-1)}^{n-1}\sup_t|\mathrm{cov}(X_{t,n},X_{t+k,n})|,$$

*which is uniformly bounded under Assumption* 2.1.

PROOF. We have, for $x \in \mathbf{C}^n$ with $\|x\|_2 = 1$ and $\sigma_{j,k} = \Sigma_{njk}$,

$$\|\Sigma_n^{1/2}x\|_2^2 = \sum_{j,k=1}^{n}\bar{x}_j\sigma_{jk}x_k \leq \Sigma_j\Sigma_k\bar{x}_j\sigma_{j,j+k}x_{j+k}.$$

An application of the Cauchy–Schwarz inequality gives the upper bound. The bound for the right-hand side follows from [14], Proposition 4.2. $\square$

We now continue with the proof of Theorem 4.1:
Let $B_n := \Sigma_n^{1/2}U_n(\frac{1}{2\pi}\phi)\Sigma_n^{1/2}$ and $\underline{Y}_n := \Sigma_n^{-1/2}\underline{X}_n \sim \mathcal{N}(0,I_n)$. We have

$$\tilde{E}_n(\phi) = n^{-1/2}\left[\underline{X}_n'U_n\left(\frac{1}{2\pi}\phi\right)\underline{X}_n - \mathrm{tr}\left\{U_n\left(\frac{1}{2\pi}\phi\right)\Sigma_n\right\}\right]$$
$$= n^{-1/2}[\underline{Y}_n'B_n\underline{Y}_n - \mathrm{tr}(B_n)].$$

Since $B_n$ is real and symmetric, there exists an orthonormal matrix $U = U_n$ with $U'U = UU' = I_n$ and $U'B_nU = diag(\lambda_{1,n},\ldots,\lambda_{n,n})$. Let $\underline{Z}_n := U'\underline{Y}_n \sim \mathcal{N}(0,I_n)$. We have

$$\tilde{E}_n(\phi) = n^{-1/2}[\underline{Z}_n'U'B_nU\underline{Z}_n - \mathrm{tr}(B_n)] = n^{-1/2}\sum_{i=1}^{n}\lambda_{i,n}(Z_i^2 - 1).$$

For $L$ and $R^2$ as defined in Proposition A.1, we obtain, with Lemma 5.2 and Lemma 5.3,

$$L = \max\{\lambda_{1,n},\ldots,\lambda_{n,n}\} = \left\|\Sigma_n^{1/2}U_n\left(\frac{1}{2\pi}\phi\right)\Sigma_n^{1/2}\right\|_{\mathrm{spec}}$$

$$\leq \left\|U_n\left(\frac{1}{2\pi}\phi\right)\right\|_{\mathrm{spec}} \|\Sigma_n^{1/2}\|_{\mathrm{spec}}^2 \leq K\rho_\infty(\phi)$$

and

$$R^2 = \frac{1}{n}\sum_{i=1}^n \lambda_{i,n}^2 = \frac{1}{n}\left\|\Sigma_n^{1/2} U_n\left(\frac{1}{2\pi}\phi\right)\Sigma_n^{1/2}\right\|_2^2$$

$$\leq \frac{1}{n}\|\Sigma_n^{1/2}\|_{\mathrm{spec}}^4 \|U_n(\phi)\|_2^2 \leq K\rho_{2,n}(\phi)^2.$$

Proposition A.1 now implies (36) and (37). Assertion (38) follows from [14], Lemma 4.3(i). The relation $\rho_{2,n}(\phi)^2 \leq \rho_2(\phi)^2 + \frac{1}{n}\rho_\infty(\phi)\tilde{v}(\phi)$ [see Remark 4.2(i)] has been proven in (65). Furthermore,

$$\begin{aligned}
\mathrm{P}(|E_n(\phi)| \geq \eta) &\leq \mathrm{P}(|\tilde{E}_n(\phi)| \geq \eta/2) + \mathrm{P}(\sqrt{n}|\mathrm{E}F_n(\phi) - F(\phi)| \geq \eta/2) \\
&\leq c_1 \exp\left(-c_2 \frac{\eta^2}{\rho_2(\phi)^2 + \frac{\eta\rho_\infty(\phi)}{\sqrt{n}} + \frac{1}{n}\rho_\infty(\phi)\tilde{v}(\phi)}\right) \\
&\quad + c_1' \exp\left(-c_2' \frac{\eta^2}{\frac{\eta}{n^{1/2}}(\rho_\infty(\phi) + \tilde{v}(\phi))}\right),
\end{aligned}$$
(66)

which implies the assertion of Remark 4.2(ii). □

PROOF OF THEOREM 4.3. We only prove the result for $d = 1$. The necessary modifications for arbitrary $d > 0$ are obvious. Let $B_n = \{\max_{t=1,\ldots,n} |X_{t,n}| \leq c\sqrt{\log n}\}$, where $c$ is the constant from Lemma A.4. This lemma says that $\lim_{n\to\infty} P(B_n) = 1$. Let $\mathcal{B}_n(\frac{\eta}{8Ln\log n})$ be the smallest approximating set at level $\frac{\eta}{8Ln\log n}$ according to the definition of the covering numbers so that $\#\mathcal{B}_n(\frac{\eta}{8Ln\log n}) = N(\frac{\eta}{8Ln\log n}, \Phi, \rho_2)$. For $\phi \in \Phi$, let $\phi^* \in \mathcal{B}_n(\frac{\eta}{8Ln\log n})$ denote the best approximation in $\mathcal{B}_n(\frac{\eta}{8Ln\log n})$ to $\phi$. With this notation, we have

$$\begin{aligned}
P\Big(\sup_{\phi\in\Phi} &|\tilde{E}_n(\phi)| > \eta; B_n\Big) \\
&\leq P\Big(\max_{\phi\in\mathcal{B}_n(\frac{\eta}{8Ln\log n})} |\tilde{E}_n(\phi)| > \eta/2\Big) \\
&\quad + P\Big(\sup_{\phi,\psi\in\Phi; \rho_2(\phi,\psi)\leq \frac{\eta}{8Ln\log n}} |\tilde{E}_n(\phi-\psi)| > \eta/2; B_n\Big) = I + II.
\end{aligned}$$
(67)

Using assumptions (41)–(43), we have

$$\begin{aligned}
I &\leq c_1 \exp\bigg\{Ak_n \log(8Ln^2\log n/\eta) - c_2 \frac{\eta^2/4}{\tau_2^2 + \frac{\eta\rho_\infty}{2\sqrt{n}} + \frac{\rho_\infty\tilde{v}}{n}}\bigg\} \\
&\leq c_1 \exp\bigg\{Ak_n \log(8Ln^2\log n/\eta) - c_2 \frac{\eta^2/4}{3\tilde{d}\tau_2^2}\bigg\} \\
&\leq c_1 \exp\bigg\{-\frac{c_2}{24\tilde{d}}\frac{\eta^2}{\tau_2^2}\bigg\}.
\end{aligned}$$
(68)

To complete the proof, we now show that for $n$ sufficiently large, we have $II = 0$ with $B_n = \{\max_{t=1,\ldots,n} |X_{t,n}| \leq c\sqrt{\log n}\}$, where $c$ is the constant from Lemma A.4 [i.e., $P(B_n) \to 1$]. In order to see this, we replace $\phi$ by

$$(69) \qquad \phi_n^*(u, \lambda) = n \int_{u-\frac{1}{n}}^{u} \phi(v, \lambda) \, dv \qquad [\text{with } \phi(v, \lambda) = 0 \text{ for } v < 0].$$

Then, on $B_n$, we have, by using Lemma A.3, the facts that $\rho_\infty(\phi - \psi) \leq 2\rho_\infty$ and $\tilde{v}(\phi - \psi) \leq 2\tilde{v}$, as well as the definitions of $B_n$ and $L$, that

$$(70) \quad \begin{aligned} |\tilde{E}_n(\phi - \psi)| &\leq \sqrt{n} |F_n(\phi - \psi) - F_n(\phi_n^* - \psi_n^*)| \\ &\quad + \sqrt{n} |F_n(\phi_n^* - \psi_n^*) - \mathrm{E} F_n(\phi_n^* - \psi_n^*)| \\ &\quad + \sqrt{n} |\mathrm{E} F_n(\phi_n^* - \psi_n^*) - \mathrm{E} F_n(\phi - \psi)| \\ &\leq 2L v_\Sigma \frac{\log n}{\sqrt{n}} + L \rho_2(\phi - \psi) n \log n + \frac{2L}{\sqrt{n}} \tilde{v} \leq \frac{\eta}{8} + \frac{\eta}{4} + \frac{\eta}{8} = \frac{\eta}{2}. \end{aligned}$$

For the last inequality to hold, we need $\eta \geq 2^4 L v_\Sigma \frac{\log n}{\sqrt{n}}$, which follows from (42). Hence, we have $II = 0$. $\square$

PROOF OF THEOREM 4.4. We use the quantities $B_n, \phi_n^*$ introduced in the proof of Theorem 4.3. Also, recall the definition of $L$ given in (39). Let

$$(71) \qquad \tilde{E}_n^*(\phi) = \sqrt{n}(F_n(\phi_n^*) - E F_n(\phi_n^*)).$$

On $B_n$, we have, by using Lemma A.3, that

$$|(\tilde{E}_n^* - \tilde{E}_n)(\phi)| \leq \sqrt{n} |F_n(\phi_n^*) - F_n(\phi)| + \sqrt{n} |(\mathrm{E} F_n(\phi_n^*) - \mathrm{E} F_n(\phi)|$$

$$\leq L \frac{\log n}{\sqrt{n}} v_\Sigma(\phi) + \frac{L}{\sqrt{n}} \tilde{v}(\phi) \leq \frac{\eta}{4} + \frac{\eta}{4} = \frac{\eta}{2},$$

where the last inequality follows from assumption (46). Hence,

$$P\Big(\sup_{\phi \in \Phi} |\tilde{E}_n(\phi)| > \eta, B_n\Big) \leq P\Big(\sup_{\phi \in \Phi} |\tilde{E}_n^*(\phi)| > \eta/2, B_n\Big).$$

We now prove the asserted maximal inequality for $\tilde{E}_n^*$. The general idea is to utilize the chaining device, as in [1].

First, we consider the case $\alpha \geq \frac{\eta}{8Ln\log n}$. In this case, choose $\delta_0 = \alpha$ and let $c_2 > 0$ be the constant from (37). Then there exist numbers $0 < \delta_j, j = 1, \ldots, K \leq \infty$, with $\alpha = \delta_0 \geq \delta_1 \geq \cdots \geq \delta_K = \frac{\eta}{8Ln\log n}$, such that with $\eta_{j+1} = \frac{3}{c_2} \delta_{j+1} \tilde{H}_\Phi(\delta_{j+1}), j = 1, \ldots, K$, we have

$$(72) \qquad \frac{\eta}{8} \geq \frac{24}{c_2} \int_{\frac{\eta}{2^5 Ln \log n}}^{\alpha} \tilde{H}_\Phi(s) \, ds \geq \sum_{j=0}^{K-1} \eta_{j+1}.$$

The first inequality follows from assumption (48) and the second follows by using the property $\delta_{j+1} \leq \delta_j/2$ (see below for the construction of the $\delta_j$).

For each of the numbers $\delta_j$, choose a finite subset $A_j$ corresponding to the definition of covering numbers $N(\delta_j, \Phi, \rho_2)$. In other words, the set $A_j$ consists of smallest possible number $N_j = N(\delta_j, \Phi, \rho_2)$ of midpoints of $\rho_2$-balls of radius $\delta_j$ such that the corresponding balls cover $\Phi$. Now, telescope

$$(73) \qquad \tilde{E}_n^*(\phi) = \tilde{E}_n^*(\phi_0) + \sum_{j=0}^{K-1} \tilde{E}_n^*(\phi_{j+1} - \phi_j) + \tilde{E}_n^*(\phi - \phi_K),$$

where the $\phi_j$ are the approximating functions to $\phi$ from $A_j$, that is, $\rho_2(\phi, \phi_j) < \delta_j$. Now take absolute value signs on both sides of (73), apply the triangle inequality on the right-hand side and then take the suprema on both sides. This leads to

$$P\Big(\sup_{\phi \in \Phi} |\tilde{E}_n^*(\phi)| > \eta/2, B_n\Big)$$
$$\leq P\Big(\sup_{\phi \in \Phi} |\tilde{E}_n^*(\phi_1)| > \eta/4\Big) + \sum_{j=0}^{K-1} N_j N_{j+1} \sup_{\phi \in \Phi} P(|\tilde{E}_n^*(\phi_{j+1} - \phi_j)| > \eta_{j+1})$$
$$+ P\Big(\sup_{\phi \in \Phi} |\tilde{E}_n^*(\phi - \phi_K)| > \eta/8, B_n\Big)$$
$$= I + II + III.$$

Note that the first two terms only depend on the approximating functions, and for every fixed $j$, those are finite in number. In contrast to that, the third term generally depends on infinitely many $\phi$ and, hence, this term is crucial. The way to treat it actually differs from case to case.

Hence, using the exponential inequality (80), we have, by definition of $\alpha$, that

$$(74) \qquad I \leq c_1 \exp\Big\{\tilde{H}(\alpha) - c_2 \frac{\eta}{4\tau_2}\Big\} = c_1 \exp\Big\{-\frac{c_2}{8} \frac{\eta}{\tau_2}\Big\}.$$

In order to estimate $II$, we need the exact definition of the $\delta_j$. Using the approach of Alexander [1], an appropriate choice [satisfying (72)] is

$$\delta_{j+1} = \frac{\eta}{8Ln \log n} \vee \sup\{x \leq \delta_j/2; \tilde{H}(x) \geq 2\tilde{H}(\delta_j)\}$$

and $K = \min\{j : \delta_j = \frac{\eta}{8Ln \log n}\}$. With these choices, we obtain

$$II \leq \sum_{j=0}^{K-1} c_1 \exp\Big\{2\tilde{H}(\delta_{j+1}) - c_2 \frac{\frac{3}{c_2}\delta_{j+1}\tilde{H}(\delta_{j+1})}{\delta_{j+1}}\Big\} \leq \sum_{j=0}^{K-1} c_1 \exp\{-\tilde{H}(\delta_{j+1})\}$$
$$\leq \sum_{j=0}^{K-1} c_1 \exp\{-2^{j+1}\tilde{H}(\alpha)\} = \sum_{j=0}^{K-1} c_1 \exp\Big\{-2^j \frac{c_2}{8} \frac{\eta}{\tau_2}\Big\} \leq 2c_1 \exp\Big\{-\frac{c_2}{8} \frac{\eta}{\tau_2}\Big\},$$

where the last inequality holds for $\eta \geq \frac{\tau_2}{c_2} \log 2$.

The proof of the fact that $III = 0$ is similar to the proof of $II = 0$ in Theorem 4.3. Here, we again have to use assumption (46). We omit the details.

It remains to consider the case $\alpha < \frac{\eta}{8Ln\log n}$. Here, we choose $\delta_0 = \frac{\eta}{8Ln\log n}$ and $K = 0$. Hence, $II = 0$ and we only have to deal with $I$ and $III$. Since $\tilde{H}(\delta_0) < \tilde{H}(\alpha)$, we immediately get [cf. (74)] that $I \leq c_1 \exp\{-\frac{c_2}{8}\frac{\eta}{\tau_2}\}$. The fact that $III = 0$ follows similarly to (70). $\square$

## APPENDIX: AUXILIARY RESULTS

First we prove a Bernstein inequality for $\chi_1^2$-variables which is the basis for the Bernstein inequality derived in Section 4.

PROPOSITION A.1. *Let $Z_1, \ldots, Z_n$ be independent standard normally distributed random variables and $\lambda_1, \ldots, \lambda_n$ be positive numbers. Define*

$$(75) \qquad R^2 = \frac{1}{n}\sum_{i=1}^n \lambda_i^2 \quad and \quad L = \max\{\lambda_1, \ldots, \lambda_n\}.$$

*Then we have for all $\eta > 0$,*

$$(76) \qquad P\left(\left|n^{-1/2}\sum_{i=1}^n \lambda_i(Z_i^2 - 1)\right| \geq \eta\right) \leq 2\exp\left(-\frac{1}{8}\frac{\eta^2}{R^2 + \frac{L\eta}{\sqrt{n}}}\right)$$

*and*

$$(77) \qquad P\left(\left|n^{-1/2}\sum_{i=1}^n \lambda_i(Z_i^2 - 1)\right| \geq \eta\right) \leq 6\exp\left(-\frac{1}{16}\frac{\eta}{R}\right).$$

PROOF. One possibility is a direct proof via moment generating functions. Instead, we apply a general Bernstein inequality for independent variables. It can be shown that $E|Z_i^2 - 1|^m \leq 4^{m-1}(m-1)!$. Therefore, we have for $m \geq 2$,

$$(78) \qquad \frac{1}{n}\sum_{i=1}^n \lambda_i^m E|Z_i^2 - 1|^m \leq \frac{m!}{2}(4L)^{m-2}(2R)^2.$$

For example, Lemma 8.6 in [22] now implies (76). Since $L \leq Rn^{1/2}$, (76) implies (77) [consider the cases $\eta \leq R$ and $\eta > R$ separately and keep in mind that $\exp(-x^2) \leq \exp(-x+1)$]. $\square$

Recall now the definition of $R_{\log}(g)$ given in (19).

LEMMA A.2. *Let $\mathcal{F}_n$ be such that Assumption 2.4(c) holds. Then we have*

$$\sup_{g \in \mathcal{F}_n} |R_{\log}(g)| = O\left(\frac{v_\Sigma}{M_* n}\right).$$

PROOF. We have, with $\phi = 1/g$,

$$|R_{\log}(g)| = \left|\frac{1}{4\pi} \int_{-\pi}^{\pi} \left[\frac{1}{n}\sum_{t=1}^{n} \log \phi\left(\frac{t}{n}, \lambda\right) - \int_0^1 \log \phi(u, \lambda)\, du\right] d\lambda\right|$$

$$\leq \frac{1}{4\pi} \int_{-\pi}^{\pi} \sum_{t=1}^{n} \int_0^{1/n} \left|\log \phi\left(\frac{t}{n}, \lambda\right) - \log \phi\left(\frac{t-1}{n} + x, \lambda\right)\right| dx\, d\lambda$$

$$\leq \frac{1}{4\pi} \int_{-\pi}^{\pi} \sum_{t=1}^{n} \int_0^{1/n} \sup_{u,\nu} |\phi(u,\lambda)^{-1}| \left|\phi\left(\frac{t}{n}, \lambda\right) - \phi\left(\frac{t-1}{n} + x, \lambda\right)\right| dx\, d\lambda$$

$$= O(n^{-1}) \frac{v_\Sigma(\phi)}{M_*}. \qquad \square$$

In the proof of Theorem 4.4, we used $\tilde{E}_n(\phi_n^*)$ instead of $\tilde{E}_n(\phi)$, where

$$\phi_n^*(u, \lambda) = n \int_{u - \frac{1}{n}}^{u} \phi(v, \lambda)\, dv \qquad [\text{with } \phi(v,\lambda) = 0 \text{ for } v < 0].$$

The reason for doing so is that otherwise, we would have needed the exponential inequality (37) to hold with $\rho_2(\phi)$ instead of $\rho_{2,n}(\phi)$. Such an inequality does not hold. Instead, we exploit the following property of $\phi_n^*$:

$$\begin{aligned}
\rho_{2,n}(\phi_n^*)^2 &= \frac{1}{n} \sum_{t=1}^{n} \int_{-\pi}^{\pi} \phi_n^*\left(\frac{t}{n}, \lambda\right)^2 d\lambda = \frac{1}{n}\sum_{t=1}^{n} \int_{-\pi}^{\pi} \left(n \int_{\frac{t-1}{n}}^{\frac{t}{n}} \phi(u,\lambda)\, du\right)^2 d\lambda \\
&\leq \sum_{t=1}^{n} \int_{-\pi}^{\pi} \int_{\frac{t-1}{n}}^{\frac{t}{n}} \phi^2(u, \lambda)\, du\, d\lambda = \rho_2(\phi)^2.
\end{aligned} \tag{79}$$

Since the assertion and the proof of Theorem 4.1 are for $n$ fixed, we obtain from (37)

$$P(|\tilde{E}_n(\phi_n^*)| \geq \eta) \leq c_1 \exp\left(-c_2 \frac{\eta}{\rho_2(\phi)}\right). \tag{80}$$

We note that $\rho_\infty(\phi_n^*) \leq \rho_\infty(\phi)$ and $\tilde{v}(\phi_n^*) \leq \tilde{v}(\phi)$, which is straightforward. The following properties are used in the proofs above:

LEMMA A.3. *Let $X_{t,n}$ be a Gaussian locally stationary process. Then we have, with $X_{(n)} := \max_{t=1,\ldots,n} |X_{t,n}|$,*

$$|F_n(\phi) - F_n(\phi_n^*)| \leq \frac{K_1}{n} X_{(n)}^2 v_\Sigma(\phi), \tag{81}$$

(82) $$|\mathrm{E}F_n(\phi) - \mathrm{E}F_n(\phi_n^*)| \leq \frac{K_2}{n}\tilde{v}(\phi),$$

(83) $$|F_n(\phi_n^*) - F(\phi)| \leq K_3(\sqrt{n}X_{(n)}^2 + 1)\rho_2(\phi) \quad and$$

(84) $$|F_n(\phi_n^*) - \mathrm{E}F_n(\phi_n^*)| \leq K_3(\sqrt{n}X_{(n)}^2 + 1)\rho_2(\phi).$$

PROOF. We have

$$|F_n(\phi) - F_n(\phi_n^*)| = \left|\frac{1}{n}\sum_{t=1}^n \int_{-\pi}^\pi \left(\phi\left(\frac{t}{n},\lambda\right) - \phi_n^*\left(\frac{t}{n},\lambda\right)\right)J_n\left(\frac{t}{n},\lambda\right)d\lambda\right|$$

$$\leq O(X_{(n)}^2)\sum_{t=1}^n \frac{1}{n}\sum_{k=-\infty}^\infty n\int_{\frac{t-1}{n}}^{\frac{t}{n}}\left|\widehat{\phi}\left(\frac{t}{n},-k\right) - \widehat{\phi}(u,-k)\right|du$$

$$\leq K_1 X_{(n)}^2 \frac{1}{n}v_\Sigma(\phi).$$

Inequality (82) can be seen as follows:

$$|\mathrm{E}F_n(\phi) - \mathrm{E}F_n(\phi_n^*)| \leq \frac{1}{2\pi}\left|\frac{1}{n}\sum_{t=1}^n\sum_k n\int_{\frac{t-1}{n}}^{\frac{t}{n}}\left[\widehat{\phi}\left(\frac{t}{n},-k\right) - \widehat{\phi}(u,-k)\right]du\right.$$

$$\left.\times \mathrm{cov}(X_{[t+1/2+k/2],n},X_{[t+1/2-k/2],n})\right|.$$

Proposition 4.2 of Dahlhaus and Polonik [14] implies $\sup_t |\mathrm{cov}(X_{t,n}, X_{t+k,n})| \leq \frac{K}{\ell(k)}$, which means that the above expression is bounded by $\frac{K\tilde{v}(\phi)}{n}\sum_k \frac{1}{\ell(k)} \leq \frac{K_2}{n}\tilde{v}(\phi)$ for some $K > 0$ independent of $n$. For the proof of (83) and (84), we estimate all terms separately by using the Cauchy–Schwarz inequality and the Parseval equality:

$$|F_n(\phi_n^*)| \leq \frac{1}{n}\sum_{t=1}^n \int_{-\pi}^\pi \left|\phi_n^*\left(\frac{t}{n},\lambda\right)J_n\left(\frac{t}{n},\lambda\right)\right|d\lambda$$

$$\leq K\rho_{2,n}(\phi_n^*)\left(\frac{1}{n}\sum_{t=1}^n \sum_{k:\,1\leq[t+1/2\pm k/2]\leq n}[X_{[t+1/2+k/2],n}X_{[t+1/2-k/2],n}]^2\right)^{1/2}$$

$$\leq K_3\rho_2(\phi)\sqrt{n}X_{(n)}^2.$$

Similarly, we obtain $|F(\phi)| \leq K_3\rho_2(\phi)$ and $|\mathrm{E}F_n(\phi_n^*)| \leq K_3\rho_2(\phi)$. □

LEMMA A.4. *Let $X_{t,n}$ be a Gaussian locally stationary process. Then there exists a $c > 0$ such that*

$$\mathbf{P}\Big(\max_{t=1,\ldots,n}|X_{t,n}| \geq c\sqrt{\log n}\Big) \to 0.$$

PROOF. We have, for some $v^* \in \mathbf{R}$, $v_{t,n} := \text{var}(X_{t,n}) = \sum_{j=-\infty}^{\infty} a_{t,n}(j)^2 \leq v^*$ uniformly in $t$ and $n$. Since $X_{t,n}$ is Gaussian, this implies for $\tilde{c} > \sqrt{2}$,

$$\mathbf{P}\left(\max_{t=1,\ldots,n}\left|\frac{X_{t,n}}{\sqrt{v^*}}\right| \geq \tilde{c}\sqrt{\log n}\right) \leq \mathbf{P}\left(\max_{t=1,\ldots,n}\left|\frac{X_{t,n}}{\sqrt{v_{t,n}}}\right| \geq \tilde{c}\sqrt{\log n}\right)$$
$$\leq \sum_{t=1}^{n} \exp\left(-\frac{\tilde{c}^2 \log n}{2}\right) \leq n^{1-\tilde{c}^2/2} \to 0. \quad \square$$

**Acknowledgments.** We are grateful to Sebastian van Bellegem for helpful comments and to the Associate Editor and two anonymous referees for careful reading of the manuscript and for their constructive criticism which led to significant improvements.

INSTITUT FÜR ANGEWANDTE MATHEMATIK
UNIVERSITÄT HEIDELBERG
IM NEUENHEIMER FELD 294
69120 HEIDELBERG
GERMANY
E-MAIL: dahlhaus@statlab.uni-heidelberg.de

DEPARTMENT OF STATISTICS
UNIVERSITY OF CALIFORNIA
DAVIS, CALIFORNIA 95616-8705
USA
E-MAIL: polonik@wald.ucdavis.edu